\documentclass[letterpaper, twoside, 11pt]{amsart}
\usepackage[T2A]{fontenc}
\usepackage[utf8]{inputenc}	
\usepackage[main=english, russian]{babel}
\usepackage{geometry}
\usepackage{comment}

\usepackage{xcolor}
\definecolor{forestgreen}{HTML}{228B22}
\definecolor{royalblue}{HTML}{4169E1}

\usepackage[unicode, pdftex]{hyperref}
\hypersetup{
	colorlinks   = true, 
	urlcolor     = blue, 
	linkcolor    = royalblue, 
	citecolor    = forestgreen 
}

\usepackage[noabbrev,capitalise,nameinlink]{cleveref}

\usepackage[abbrev,alphabetic,backrefs,initials,nobysame]{amsrefs}

\usepackage{amssymb}
\usepackage{tikz-cd}
\usepackage[all]{xy}
\usepackage{mathrsfs}
\usepackage{xfrac}

\newtheorem{theorem}{Theorem}[section]
\newtheorem{proposition}[theorem]{Proposition}
\newtheorem{lemma}[theorem]{Lemma}

\newtheorem{conjecture}[theorem]{Conjecture}

\theoremstyle{definition}
\newtheorem{definition}[theorem]{Definition}

\newtheorem{example}[theorem]{Example}

\theoremstyle{remark}
\newtheorem{remark}[theorem]{Remark}

\newtheorem*{example*}{Example}
\newtheorem*{question*}{Question}

\theoremstyle{plain}
\newcounter{intro}

\newtheorem{intro-conjecture}[intro]{Conjecture}
\newtheorem{intro-corollary}[intro]{Corollary}
\newtheorem{intro-theorem}[intro]{Theorem}
\newtheorem*{intro-corollary*}{Corollary}
\newtheorem*{intro-conjecture*}{Conjecture}

\newcommand{\Q}{\mathbb{Q}}
\newcommand{\R}{\mathbb{R}}
\newcommand{\Ko}{\mathbb{C}}
\newcommand{\Zi}{\mathbb{Z}}

\newcommand{\Cc}{\mathfrak{C}}
\newcommand{\Dd}{\mathfrak{D}}

\newcommand{\X}{\mathcal{X}}

\newcommand{\B}{\mathcal{B}}

\newcommand{\A}{\mathcal{A}}

\newcommand{\Y}{\mathcal{Y}}

\newcommand{\Os}{\mathcal{O}}

\newcommand{\Pp}{\mathbb{P}}
\newcommand{\al}{\alpha}
\newcommand{\Aut}{\textnormal{Aut}}

\newcommand{\supp}{\textnormal{supp}\,}

\newcommand{\Iso}{\textnormal{Iso}}
\newcommand{\Wlc}{\textnormal{Wlc}}
\newcommand{\wlc}{\textnormal{wlc}}

\let\emptyset\varnothing

\title{Constructible properties of minimal models}
\author[D. V. Serebrennikov]{Daniil Serebrennikov}
\address{Department of Mathematics, Johns Hopkins University, Baltimore, MD 21218, USA.}
\email{dserebr1@jhu.edu}

\begin{document}
	\begin{abstract}
		We prove that, given a projective minimal model and a projective family of good minimal models, the subset of the base parametrizing fibers isomorphic (or birationally equivalent) to the given minimal model is constructible in the Zariski topology. Together with recent boundedness results for polarized log Calabi--Yau fibrations, this constructibility theorem provides new finiteness results predicted by the relative Kawamata--Morrison--Totaro cone conjecture.
	\end{abstract}
	\maketitle
	\tableofcontents
	\section{Introduction}
		\noindent Throughout, we work over an algebraically closed field $k$ of characteristic 0.
		
		Let $f: \X\to S$ be a projective morphism between quasi-projective varieties, and $\mathbf{P}$ a property of schemes over $k$. It is natural to ask the following question.
		\begin{question*}
			 Is the set $S_{\mathbf{P}} := \{s\in S: \mathbf{P}(\X_s) \text{ holds}\}$ constructible?
		\end{question*}
		
		\noindent Recall that a subset\footnote{In the introduction, we consider only closed points of the base.} of $S$ is \textit{constructible} if it is a finite union of locally closed subsets in the Zariski topology. For instance, Grothendieck \cite{EGAIV} proved that many basic geometric properties, such as smoothness, normality, and irreducibility, define constructible loci (see \cite[Appendix E]{GW20} for a summary of such results). These results naturally lead to the analogous question for the properties of being isomorphic or birationally equivalent to a fixed variety. For a fixed projective variety $X$, set $$\Iso_X(\X/S) := \{s\in S: \X_s\approx X\}, \ \text{Bir}_X(\X/S) := \{s\in S: \X_s\sim_{\text{bir}} X\}.$$
		\begin{question*}
			Are the sets $\text{Iso}_X(\X/S)$ and $\text{Bir}_X(\X/S)$ constructible?
		\end{question*}
		\noindent In general, the sets $\text{Iso}_X(\X/S)$ and $\text{Bir}_X(\X/S)$ are at most countable unions of constructible subsets, cf. \cite{FGA05, BBG16}. For example, in \cite{HPT18} the authors constructed a smooth projective family $\X/S$ of  complex 4-folds such that the set $\text{Bir}_{\Pp_\Ko^4}(\X/S)$ is a countable union of closed subsets, which is dense in $S$ even for the Euclidean topology. Nevertheless, in our previous work \cite{Ser26b}, we proved the constructibility of both loci for smooth projective  families $\X/S$ of K-torsion varieties over $\Ko$, i.e. $\omega_{\X_s}$ is torsion for every $s\in S$, via the Hodge theory \cite{BBT23, BFMT25}. Eventually, constructibility for $\text{Iso}_X(\X/S)$ fails  for non-minimal varieties in general as the following example shows \cite[Section 1.6]{Huy18}.
		\begin{example*}
			Suppose $k = \Ko$. Let $S$ be a projective smooth surface such that $K_S$ is nef, $K_S^2 = 0$, and $\text{Aut}(S)$ is a countable group (see \cite{SS19}). Let $\Delta\subseteq S\times_k S$ be the diagonal, and $\X = \text{Bl}_\Delta(S\times_k S)$ be its blow-up. We claim that $\text{Iso}_{\X_0}(\X/S) = \text{Aut}(S)\cdot o$ for any fixed $o\in S$. For every $s\in \text{Iso}_{\X_0}(\X/S)$ an isomorphism $\X_o\to \X_s$ induces a birational map $S\dashrightarrow S$ which is an isomorphism because $S$ is a minimal surface. Hence, $\text{Iso}_{\X_0}(\X/S) \subseteq \text{Aut}(S)\cdot o$. Conversely, every such automorphism lifts to an isomorphism of the two blow-ups. Now, let $\text{Stab}(\varphi) = \{ s\in S: \varphi(s) = s\}$, and choose $o\in \bigcap_{\varphi\in \text{Aut}(X)}S\smallsetminus \text{Stab}(\varphi)$. Then the orbit $\text{Aut}(S)\cdot o = \text{Iso}_{\X_0}(\X/S)$ is a countable set that is not constructible.
		\end{example*}
		
		Interestingly, we showed \cite[Section 5.3]{Ser26b} the constructibility of  $\text{Iso}_X(\X/S)$ is also sensitive to singularities of $X$. From  the Minimal Model Program point of view \cite[Section 2.1]{Kol13}, it is natural to study singularities in the framework of log pairs $(X,B)$, here $X$ is a normal variety and $B$ is an $\R$-divisor such that $K_X + B$ is $\R$-Cartier. We say  a log pair $(X, B)$ is a \textit{0-pair} if $X$ is proper and $K_X + B\equiv 0$. More generally, we proved that the loci $\text{Iso}_{(X, B)}(\X/S, \B)$ and $\text{Bir}_{(X, B)}(\X/S, \B)$ are constructible for projective elementary families $(\X/S, \B)$ of 0-pairs $(\X_s, \B_s)$ with at worst Kawamata log terminal singularities (klt) and $\Q$-coefficients. These constructibility results also imply finiteness statements for $0$-pairs \cite[Section 5]{Ser26a}. For instance, on a fixed projective klt 0-pair $(X,B)$ with $\Q$-coefficients, there are at most finitely many numerical classes of ample Cartier divisors of fixed volume up to the action of $\text{Aut}(X, B)$.
		
		These finiteness results \cite[Appendix A]{Ser26a} are predicted by the Kawamata--Morrison--Totaro cone conjecture \cite{Mor93, Kaw97, Tot10}. Surprisingly, even the constructibility itself provides a new evidence for the cone conjecture as the latter implies the former \cite[Corollary]{Ser26b}. We also provided examples of lc 0-pairs $(X, B)$, including plt examples, for which the analogous cone conjecture fails and the locus $\Iso_{(X,B)}(\mathcal{X}/S,\mathcal{B})$ is not constructible in suitable projective families \cite[Examples~5.11 and~5.13]{Ser26b}. The reason of failure is that these pairs admit infinitely many numerical classes of ample divisors of fixed volume, even up to $\text{Aut}(X)$.  In this paper, we extend our previous results to all projective \textit{good weak log canonical (wlc) models } $(X,B)$ with at worst klt singularities, that is, $B$ is an effective $\R$-divisor such that $K_X + B$ is nef and semiample.
		
		\begin{intro-theorem}
			\label{intro-theorem: constructibility}
			Let $(\X/S, \B)$ be a projective elementary family of good klt wlc models, and $(X, B)$ a good klt wlc model. Then the following subsets are constructible:
			\begin{equation*}
				\Iso_{(X, B)}(\X/S, \B), \textnormal{Bir}_{(X, B)}(\X/S, \B)\subseteq S.
			\end{equation*}
		\end{intro-theorem}
		
		Our previous work \cite{Ser26b} relied on transcendental methods over $\Ko$. In contrast, the arguments in the present paper are entirely algebraic. Nevertheless, we have to justify the strong Lefschetz principle. This justification is necessary because the constructibility of a locus defined by a geometric property can depend on the ground field, as the following example illustrates.

		\begin{example*}
			Let $\X/S$ be the Legendre family of smooth projective genus 1 curves over $k\in \{\overline{\Q}, \Ko\}$, here $S = \mathbb{A}^1_k\smallsetminus\{0, 1\}$. We say a property $\mathbf{P} = \mathbf{P}(Z)$ holds for $Z$ if there is a finite map $Z\to \Pp^1_k$ ramified at three points only. According to Belyi's theorem, we have $S_\mathbf{P} = \{s\in S: \X_s \text{ is defined over } \overline{\Q}\}$. Then $S_{\mathbf{P}} = S$ if $k = \overline{\Q}$. In contrast, if $k = \Ko$ the locus $S_\mathbf{P}\subset S$ is countable that is not constructible.
		\end{example*}
		
		For applications, we provide further evidences (finiteness of big and nef/movable orbits) for the  relative Kawamata--Morrison--Totaro cone conjecture as well as for the conjecture concerning finiteness of minimal models. Recall that wlc models $(X, B)$ and $(X', B')$ are \textit{0-equivalent} if $X \sim_\text{bir}~X'$, and the pullbacks of their log canonical divisors $K_X+ B$ and $K_{X'} + B'$ coincide on any higher birational model.
		
		\begin{intro-conjecture*}[{\cite[Conjecture A.2]{Ser26b}}]
			\label{intro-conj: KWS}
			The number of log isomorphism classes of projective klt wlc models in a fixed 0-equivalence class is finite.
		\end{intro-conjecture*}
		
		For terminal $\Q$-factorial minimal models, this conjecture is due to Kawamata and Matsuki \cite[Conjecture~12.3.6]{Mat02}. Its generalization to log pairs was proposed and proven for log surfaces in \cite{Ser26a}, though this case follows from \cite{Ste76, Tot10}.  In fact, the conjecture has many applications in birational geometry, and it is known for wlc models of general type \cite{BCHM10}. In this paper, we provide a partial resolution of the conjecture above by reducing it to existence of very ample Cartier divisors of bounded degree on each wlc model (see \cref{def: bounded_polarization}).
		
		\begin{intro-theorem}
			\label{intro-theorem: finiteness_and_boundedness}
			Let $\Cc$ be a collection of 0-equivalent projective good klt wlc models. Then the set of log isomorphism classes in $\Cc$  is finite if and only if $\Cc$ admits a bounded polarization.
		\end{intro-theorem}
		
		For a proper surjective morphism $X/Y$ between varieties, let $\text{NS}(X/Y)$ denote the group of Cartier divisors on $X$ modulo the numerical equivalence over $Y$, and $X_\eta$ be the generic fiber. Then the group $\text{Aut}(X/Y)$ has a natural right action on $\text{NS}(X/Y)$ defined by pullbacks. Moreover, $\Aut(X/Y)$ also naturally acts on $\text{NS}(X_\eta)$ using the restriction map $\text{NS}(X/Y)  \to \text{NS}(X_\eta)$. By $[A_\eta]$ we denote the restriction of $[A]\in \text{NS}(X/Y)$ to $\text{NS}(X_\eta)$.

		\begin{intro-theorem}
			\label{intro-theorem: finiteness_big_nef}
			Let $(X, B)$ be a projective good klt wlc model, $Y$ its log canonical model, and $X_\eta$ the generic fiber of $X/Y$.  Then for every integer $v\in \mathbb{Z}$, the number of orbits for the natural action of $\textnormal{Aut}(X/Y, B)$ on the following set is finite:
			\begin{equation*}
				\{ \ [A_\eta]\in \textnormal{NS}(X_\eta): \textnormal{vol}(A_\eta) = v,\  A \textnormal{ is a Cartier divisor on $X$ which is big and nef  over $Y$ } \}.
			\end{equation*}
		\end{intro-theorem}
		
		Similarly, 	for a proper surjective morphism $X/Y$ between varieties by $\text{NS}_W(X/Y)$ we denote the group of Weil divisors on $X$ modulo the numerical equivalence over $Y$, which is defined when $X$ is normal and $\Q$-factorial. Let $X_\eta$ be the generic fiber. Then the group $\text{PsAut}(X/Y)$ has a natural right action on $\text{NS}_W(X/Y)$ defined via inverse strict transforms.  Moreover, $\text{PsAut}(X/Y)$ also naturally acts on $\text{NS}_W(X_\eta)$ using the restriction map $\text{NS}_W(X/Y)  \twoheadrightarrow \text{NS}_W(X_\eta)$. By $[M_\eta]$ we denote the restriction of $[M]\in \text{NS}_W(X/Y)$ to $\text{NS}_W(X_\eta)$, see \cref{sec: applications}.

		\begin{intro-theorem}
			\label{intro-theorem: finiteness_big_movable}
			Let $(X, B)$ be a projective good klt wlc model, $Y$ its log canonical model, and $X_\eta$ the generic fiber of $X/Y$.  In addition, we suppose $X$ is $\Q$-factorial. Then for every number $v\in \mathbb{Q}$, the number of orbits for the natural action of $\textnormal{PsAut}(X/Y, B)$ on the following set is finite:
			\begin{equation*}
				\{ \ [M_\eta]\in \textnormal{NS}_W(X_\eta): \textnormal{vol}(M_\eta) = v,\  M \textnormal{ is a Weil divisor on $X$ which big and movable  over $Y$ } \}.
			\end{equation*}
		\end{intro-theorem}

	\addtocontents{toc}{\protect\setcounter{tocdepth}{1}}
	\subsection*{Outline of proofs} 
	\begin{proof}[Sketch of the proof for \cref{intro-theorem: constructibility}]
		Proceeding by Noetherian induction on closed subsets, we can suppose $\text{Bir}_{(X, B)}(\X/S, \B)$ is dense in $S$. After shrinking the base, we construct the log canonical model $\Y/S$ for $(\X/S, \B)$. Since the log canonical model is unique for given 0-equivalent pairs and the subset of the base corresponding to isomorphic \text{polarized pairs} is  constructible, we reduce the problem to the case when $\Y = Y\times_k S$ where $Y$ is the log canonical model for $(X, B)$. Working over the geometric generic point $K$ of $Y$, we obtain an isotrivial lc-trivial fibration over $S_K$ via \cite{Ser26b} and the strong Lefschetz principle. From descent theory it follows that after an \'etale base change the family $(\X/S, \B)$ is crepant birationally equivalent to a trivial one. From rigidity results  \cite{Kol26} for small maps and birational contractions it follows that $(\X/S, \B)$ is trivial as well, as required.
	\end{proof}
	\begin{proof}[Sketch of the proof for \cref{intro-theorem: finiteness_and_boundedness}]
		We reduce the problem to the case when the ground field $k$ is uncountable via the strong Lefschetz principle. By the standard Hilbert--Chow schemes argument together with Noetherian induction, it is suffices  to consider a projective family $(\X/S, \B)$ of good klt wlc models such that $\text{Bir}_{(X, B)}(\X/S, B)$ is dense in $S$. By \cref{intro-theorem: constructibility}, we can assume $S = \text{Bir}_{(X, B)}(\X/S, B) = \bigsqcup_\al \text{Iso}_{(X_\al, B_\al)}(\X/S, \B)$  where the union runs over the distinct log isomorphism classes occurring in the family. In fact, this decomposition is at most countable. Since every set in this decomposition is constructible, it follows that one subset $\text{Iso}_{(X,B)}(\X/S, \B)$ includes a dense open subset of $S$ because $k$ is uncountable, as needed.
	\end{proof}  
	
	Using Shokurov's log geography together with \cref{intro-theorem: finiteness_and_boundedness}, we can reduce both  \cref{intro-theorem: finiteness_big_nef} and \cref{intro-theorem: finiteness_big_movable} to the problem that the number of orbits for the natural action of $\text{Aut}(X/Y, B)$ on the following set is finite:
	\begin{equation*}
		\{[A_\eta]\in \textnormal{NS}(X_\eta): \textnormal{vol}(A_\eta) = v,\  A \textnormal{ is a Cartier divisor on $X$ which is ample over } Y\}.
	\end{equation*}
	The main difficulty is to construct global polarizations of	uniformly bounded volume while controlling the	prescribed numerical classes on the generic fiber. The necessary techniques to address this issue were established in \cite{Bir23}. Recently, a generalization to polarized log Calabi--Yau fibrations was treated in \cite{JJZ26, Zhu26}. See also \cite{BDCS24, Fil24, FHS25, HH25, Jia25, Zhu25}.
	
	Once the required boundedness  have	been established, we apply \cref{intro-theorem: constructibility} to obtain finitely many trivial families in which the given divisors varies algebraically. Their numerical classes are constant in each family, giving the desired finiteness of orbits.

	\subsection*{Acknowledgments}
	The author is deeply grateful to his advisor, Vyacheslav Shokurov, for inspiring discussions and encouragement. The author also thanks Stefano Filipazzi for several very fruitful discussions.

	\section{Preliminaries}
	\subsection{Notation and conventions.}
	Throughout, we work in the category of Noetherian schemes over an algebraically closed field $k$ of characteristic 0, and follow the terminology introduced in \cite{Ser26b} unless stated otherwise. For reader convenience, we shall recall only main definitions from \textit{loc. cit.} and adapt some of them to the relative setting.
	
	By a variety we mean an integral separated scheme of finite type over a field.
	On a variety, a prime divisor is understood in the Weil sense, i.e. an integral closed subscheme of codimension~$1$. In every part of the paper, for an $\mathbb{R}$-divisor $D$ the sum $D=\sum_i d_i D_i$ assumes that all $D_i$ are distinct prime divisors, and $d_i\in\mathbb{R}$. Otherwise, we shall say $D=\sum_i d_i D_i$ is a \textit{decomposition}.  A contraction $f: X\to Y$ is a proper surjective morphism with connected fibers.
	
	\subsection{Pairs.} A \textit{pair} $(X,B)$ consists of a variety $X$ and an $\mathbb{R}$-divisor $B$ on $X$. We say that $(X,B)$ is a \textit{log pair} if $X$ is normal, and $K_X+B$ is an $\mathbb{R}$-Cartier $\R$-divisor. Suppose $B=\sum_i b_i B_i$, and $\Gamma\subseteq\mathbb{R}$ is a set. If $b_i\in\Gamma$ for all~$i$ then we write $\text{coeff }B\subseteq \Gamma$. Denote by $[0,1]$ the unit segment of real numbers.  The $\mathbb{R}$-divisor $B$ is called a \textit{boundary} if $\text{coeff }B\subseteq  [0,1]$.	The pairs $(X,B)$ and $(X',B')$ are said to be \textit{log isomorphic} if there exists an isomorphism $\varphi:X\to X'$ such that $\varphi_*(B)=B'$; we then write $(X,B)\approx (X',B')$. We use the notation $\cong$ instead of $\approx$ if an isomorphism is natural. We write $g: (Y, D) \to (X, B)$ if $g: Y\to X$ is a morphism such that $K_Y + D = g^*(K_X + B)$.
	
	\begin{definition}
		\label{def: sing}
		Let $(X,B)$ be a log pair with $B=\sum b_i B_i$. We say that $(X,B)$ has (at worst) the following singularities when the corresponding inequalities hold:
		\[
		\begin{gathered}
			\begin{cases}
				\text{terminal (trm)} \\
				\text{Kawamata log terminal (klt)}\\
				\text{log canonical (lc)}
			\end{cases}
			\Longleftrightarrow
			\text{dis}(X,B)
			\begin{cases}
				>0;\\
				>-1 \ \text{and}\  b_i<1 \ \forall i;\\
				\ge -1 \ \text{and}\  b_i\le 1 \ \forall i.
			\end{cases}
		\end{gathered}
		\]
	\end{definition}
	
	\begin{definition}
		\label{def: 0-pair}
		A log pair $(X,B)$ with a morphism $X/S$ is called a \textit{weakly log canonical model over $S$} if the following assumptions hold:
		\begin{itemize}
			\item $X/S$ is a proper morphism between varieties, and $B$ is a boundary.
			\item $(X,B)$ has (at worst) log canonical singularities.
			\item $K_X+B$ is nef over $S$.
		\end{itemize}
		Then we  write $(X, B)/S$ and simply say it is a wlc model. If $(X,B)$ has klt (resp. terminal) singularities, then $(X,B)/S$ is called a \textit{klt wlc model} (resp. a \textit{trm wlc model}). If, in addition, $K_X+B\equiv_S 0$, then $(X,B)/S$ is called a \textit{0-pair}. When $S$ is a point, we omit it from the notation.
	\end{definition}
	
	\begin{definition}
		A wlc model $(X, B)/S$ is called \textit{good} if $K_X + B$ is semiample over $S$.
	\end{definition}
	
	\begin{conjecture}[Abundance]
		Every projective wlc model $(X, B)$ is good.
	\end{conjecture}
	The above conjecture predicts that every wlc model is a  0-pair over its log canonical model. In practice, this allows one to reduce study of wlc models to families of  0-pairs.
	\begin{definition}
		\label{def: 0-class}
		We say that a log pair $(X_\alpha,B_\alpha)$ is \textit{crepant birationally equivalent} or \textit{0-equivalent} to $(X,B)$ if the varieties  $X_\alpha$ and $X$ are birationally equivalent, and for each common log resolution $Y$ of these pairs there is an $\R$-divisor $D$ on $Y$ such that the following equalities hold:
		\[
		\begin{gathered}
			K_Y+D=g^*(K_X+B)=g_\alpha^*(K_{X_\alpha}+B_\alpha),\\
			B=g_*D,\quad B_\alpha=(g_\alpha)_*D,
		\end{gathered}
		\]
		where $g:Y\to X$ and $g_\alpha:Y\to X_\alpha$ are the corresponding log resolutions. The class of all log pairs crepant birationally equivalent to $(X,B)$ is called the \textit{0-class} of the log pair $(X,B)$.
	\end{definition}
	
	When log pairs $(X, B)$ and $(X', B')$ are 0-equivalent, we write $(X_\al, B_\al)\sim_{\text{bir}} (X, B)$. If, in addition, the varieties $X$ and $X_\al$ are isomorphic in codimension 1, we write $(X_\al, B_\al)\sim_{\text{iso}^1} (X, B)$.
	
	\subsection*{Families}
		Let $f:\X\to S$ be a morphism of schemes and $s\in S$ a point. We denote the scheme-theoretic fiber by $\X_s=\X\times_S s$. If $\A$ is an invertible sheaf on $\X$, then its restriction to the fiber $\X_s$ is denoted by $\A_s=\A|_{\X_s}$. When no ambiguity arises, we denote a morphism $f: \X \to S$ simply by $\X/S$. If $\A$ is (very) ample over $S$, we  say that $\A$ is an (very) ample invertible sheaf on $\X/S$.
	
	\begin{definition}
		Suppose $f: \X \to S$ is a proper flat morphism between varieties, and the base $S$ is regular. We say that $\X/S$ is a \textit{family of varieties} if every fiber $\X_s$ is a variety. If, in addition, the morphism $f$ is projective (resp. smooth), then $\X/S$ is called a projective (resp. smooth) family of varieties.
	\end{definition}
	
	\begin{definition}
		\label{def: elementary_family}
		Let $(\X,\B)$ be a pair and $f:\X\to S$ be a family of varieties. Suppose $\B=\sum_{i} b_i\B_i$ is a decomposition of the $\mathbb{R}$-divisor into distinct prime divisors. We say that $(\X/S,\B)$ is an \textit{elementary family} if for every point $s\in S$ the following hold:
		\begin{enumerate}
			\setcounter{enumi}{-1}
			\item The fiber $\X_s$ is a normal variety.
			\item The restriction $f|_{\B_i}:\B_i\to S$ is flat for all $i$.
			\item The closed subscheme $(\B_i)_s=\B_i\times_S s$ is a prime divisor on $\X_s$ for all $i$.
			\item The prime divisors $(\B_i)_s$ and $(\B_j)_s$ are distinct whenever $i\ne j$.
		\end{enumerate}
		In this situation, we put $\B_s=\sum_i b_i(\B_i)_s$. If, in addition, the morphism $f$ is projective, then $(\X/S,\B)$ is called a \textit{projective elementary family}.
	\end{definition}
	
	\begin{definition}
		A \textit{polarized pair} $(X, B; A)$ is a triple such that $(X, B)$ is a pair, and $A$ is an ample invertible sheaf on $X$. We say that two polarized pairs $(X, B; A)$ and $(X^\prime, B^\prime; A^\prime)$ are isomorphic if there exists an log isomorphism $\varphi: (X, B)\to (X^\prime, B^\prime)$ such that $A \approx \varphi^*A^\prime$, the latter is an isomorphism of invertible sheaves on $X$. In this situation, we write $(X, B; A) \approx (X^\prime, B^\prime; A^\prime)$.
	\end{definition}
	
	\begin{definition}
		\label{def: bounded_pairs}
		A class $\Cc$ of polarized pairs $(X_\alpha,B_\alpha; A_\al)$ is called \textit{bounded} if there exist finitely many elementary families $(\X^{(j)}/S^{(j)},\B^{(j)})$ together with an ample invertible sheaf $\A^{(j)}$ on $\X^{(j)}/S^{(j)}$ such that:
		\begin{itemize}
			\item For every $(X_\alpha,B_\alpha; A_\al)$ in $\Cc$ there is an index $j$ and a closed point $s_\alpha\in S^{(j)}$ such that $(X_\alpha,B_\alpha; A_\al)\approx~(\X^{(j)}_{s_\alpha},\B^{(j)}_{s_\alpha}; \A^{(j)}_{s_\alpha}).$
			\item For each index $j$ and every closed point $s\in S^{(j)}$ there exists a polarized pair $(X_\alpha,B_\alpha; A_\al)$ in $\Cc$ such that $(\X^{(j)}_s,\B^{(j)}_s, \A^{(j)}_{s})\approx~(X_\alpha,B_\alpha; A_\al)$.
		\end{itemize}
	\end{definition}
	
	\begin{definition}
		\label{def: bounded_polarization}
		A class $\Cc$ of projective pairs $(X_\alpha,B_\alpha)$ of fixed dimension $d\in \mathbb{Z}_{>0}$ \textit{admits a bounded polarization} if there exists $N\in \mathbb{Z}_{> 0}$ and a finite subset $\Gamma\subset\mathbb{R}$ such that $B_\alpha\in\Gamma$, and each $X_\alpha$ admits a very ample invertible sheaf $A_\al$ on $X_\al$ satisfying the following inequalities:
		\[
		(B_\alpha)_{\mathrm{red}}\cdot A_\al^{d-1}\le N,\qquad
		A_\al^{d}\le N,
		\]
		where $(B_\alpha)_{\mathrm{red}}=\sum_{V\in\supp(B_\alpha)} V$, $\ \supp(B_\alpha)=\{\,B_{i\alpha}\mid b_{i\alpha}\ne 0\,\}$, and $B_\alpha=\sum_i b_{i\alpha}B_{i\alpha}$. The class of triples$(X_\al, B_\al; A_\al)$ is called \textit{a class of polarized pairs associated with $\Cc$}.
	\end{definition}

	\begin{proposition}[{\cite[Proposition 2.15]{Ser26b}}]
		\label{prop: bounded_subclass}
		Let $\Cc$ be a class of pairs of fixed dimension $d \in \Zi_{>0}$. Then any class of polarized pairs associated with $\Cc$ is a subclass of a bounded class of polarized pairs, whenever $\Cc$ admits a bounded polarization.
	\end{proposition}
	\begin{proof}
		The proof in \textit{loc. cit.} is valid for every algebraically closed field of characteristic 0 as the proof essentially involves only existence of the Chow and Hilbert schemes.
	\end{proof}

	\subsection{Constructibility.} To prove the main constructibility theorem in this paper, which is \cref{th: main_constr}, we need results from our previous work \cite{Ser26b} where transcendental methods were used. In contrast, current paper is based on algebraic methods only. To overcome possible subtleties, we justify the strong Lefschetz principle in this subsection directly.
	
	\begin{definition}
		Let $S$ be a Noetherian topological space. A subset $W\subseteq S$ is called \textit{constructible} if there are finitely many closed and open subsets $Z_i, U_i\subseteq S$ respectively such that $W = \bigcup_i (Z_i\cap U_i)$.
	\end{definition} 
	
	\begin{lemma}
		\label{lem: constr_fields}
		Let $\iota: S'\to S$ be a surjective morphism between Noetherian schemes. Then a subset $W\subseteq S$ is constructible if and only if $\iota^{-1}(W)$ is constructible.
	\end{lemma}
	\begin{proof}
		If $W$ is constructible, then $\iota^{-1}(W)$ is constructible because morphisms of schemes are continuous maps in the Zariski topology. Suppose $\iota^{-1}(W)$ is constructible. Note that $\iota$ is quasi-compact because $S'$ is Noetherian. Since the question is local, we can assume that $S$ is affine. Choose a finite affine open cover $S' = \bigcup_i U_i$. After the base change $\bigsqcup_i U_i \to \bigcup_i U_i$, we can assume that $S'$ if affine as well, i.e. $\iota: S'\to S$ is a  surjective spectral map between spectral spaces (see \cite[\href{https://stacks.math.columbia.edu/tag/08YF}{Tag 08YF}]{Stacks}). Now, we endow $S$ and $S'$ with the constructible topologies, that is, the coarsest topology in which every constructible subset is open and closed. It is clear that $S, S'$ are Hausdorff in that topology. Then $S\smallsetminus W = \iota(\iota^{-1}(S\smallsetminus W))$ is closed as a quasi-compact subset of the Hausdorff space, here we used the surjectivity of $\iota$. Hence, $W\subseteq S$ is open and quasi-compact subset in the constructible topology. This implies that $W$ is constructible, as required.
	\end{proof}
	Given a pair $(X, B)$ and a field extension $K/k$, by $(X, B)_K$ we denote the pair obtained by the constant field extension, that is, $(X, B)_K:=(X_K, B_K)$ where $X_K = X\times_k K$ and $B_K = B\times_k K$.
	\begin{lemma}
		\label{lem: iso_fields}
		Let $(X, B)$ and $(X', B')$ be log pairs, and $K/k$ a field extension. Then 
		\begin{enumerate}
			\item 
				$(X, B)_K \sim_{\text{bir}} (X', B')_K$ if and only if $(X, B) \sim_{\text{bir}} (X', B')$.
			\item 
				$(X, B)_K \sim_{\text{iso}^1} (X', B')_K$ if and only if $(X, B) \sim_{\text{iso}^1} (X', B')$.
			\item 
				 $(X, B)_K \approx (X', B')_K$ if and only if $X\approx X'$.
		\end{enumerate}
	\end{lemma}
	\begin{proof}
		One direction is obvious for every statement. Suppose that $(X, B)_K \sim_{\text{bir}} (X', B')_K$ over $K$, that is, there is a log pair $(Y, D)$ over $K$ with crepant contractions $f: (Y, D)\to (X, B)_K$ and $f': (Y, D)\to (X', B')_K$. Let $A_0$ be an algebra over $k$ generated by all coefficients of (finitely many) polynomials defining $(Y, D)$.  Then we can write $K = \varinjlim A$ where the direct limit runs over all finitely generated  over $k$ algebras $A$ such that $A_0\subseteq A\subseteq K$. By the standard descent argument \cite[\href{https://stacks.math.columbia.edu/tag/01ZM}{Tag 01ZM}]{Stacks}, there exists a finitely generated over $k$ algebra $A\supseteq A_0$  over which morphisms $f, f'$ descend, that is, there are morphisms $f_A: (Y, D)_A\to (X, B)_A$ and $f_A': (Y, D)_A\to (X', B')_A$. Next, we can assume that they are crepant contractions. These morphisms induce crepant contractions $f_s: (Y, D)_{\kappa(s)} \to (X, B)_{\kappa(s)}$ and $f_s': (Y, D)_{\kappa(s)}\to (X', B')_{\kappa(s)}$ for any closed point $s\in \text{Spec}\, A$. Since $k$ is algebraically closed, we obtain $\kappa(s) \cong k$. This implies $(X, B) \sim_{\text{bir}} (X', B')$ over $k$. Now, in addition, we suppose that $(X, B)_K \sim_{\text{iso}^1} (X', B')_K$ over $K$. Let $Z\subseteq X$ (resp. $Z'\subseteq X$) be the locus of indeterminacy for $(f_s)^{-1}\circ f'_s$ (resp. $(f'_s)^{-1}\circ f_s$). Since the codimension is unchanged by field extensions, we obtain that $\text{codim}_X Z = \text{codim}_{X_K} Z_K\ge 2$ and $\text{codim}_{X'} Z' = \text{codim}_{X'_K} Z_K'\ge 2$. Finally, we suppose  that there is an isomorphism  $\varphi: (X, B)_K \to (X', B')_K$ over $K$. To finish the proof, one can take $f = \text{id}_X$ and $f' = \varphi$.
	\end{proof}
	
	Let $(\X/S, \B)$ be an elementary family, and $(X,B)$ be a proper log pair. Suppose $s\in S$ is a point, $\kappa(s)$ is its residue field, and $\overline{\kappa(s)}$ is an algebraic closure for  $\kappa(s)$. By $(\X, \B)_{\overline{s}}$ we denote the geometric fiber over $s$, that is, $(\X_s, \B_s)\times_{\kappa(s)} \overline{\kappa(s)}$.  Define the following sets:
	\begin{equation*}
		\begin{gathered}
			\text{Iso}_{(X, B)}(\X/S, \B) =  \{s\in S: (\X, \B)_{\overline{s}}\approx (X, B)_{\overline{s}}\},\\
			\text{Iso}^1_{(X, B)}(\X/S, \B) =  \{s\in S: (\X, \B)_{\overline{s}}\sim_{\text{iso}^1}  (X, B)_{\overline{s}}\},\\
			\text{Wlc}_{(X, B)}(\X/S, \B) = \{s\in \text{Bir}_{(X, B)}(\X/S, \B):(\X, \B)_{\overline{s}} \text{ is a wlc model} \},\\
		\end{gathered}
	\end{equation*}
	where $\text{Bir}_{(X, B)}(\X/S, \B) = \{s\in S:(\X, \B)_{\overline{s}}\sim_{\text{bir}} (X, B)_{\overline{s}}\}$. It is clear that $$\Iso_{(X, B)}(\X/S, \B)\subseteq \Iso^1_{(X, B)}(\X/S, \B)\subseteq \text{Bir}_{(X, B)}(\X/S, \B).$$ When $(X, B)$ is a wlc model, we have one additional inclusion $ \Iso^1_{(X, B)}(\X/S, \B)\subseteq \Wlc_{(X, B)}(\X/S, \B)$.

	\begin{theorem}[{\cite{Ser26b}}]
		\label{th: 0-pairs_constr}
		Let $(\X/S, \B)$ be a projective elementary family, and $(X, B)$ a projective klt 0-pair, where $B$ is an (effective) $\Q$-divisor. Then the following subsets  are constructible: $$\Iso_{(X, B)}(\X/S, \B), \Iso^1_{(X, B)}(\X/S, \B), \Wlc_{(X, B)}(\X/S, \B)\subseteq S.$$ 
	\end{theorem}
	\begin{proof}
		Let $k_0/\Q$ be a finitely generated field over which the family $(\X/S, \B)$ and the pair $(X, B)$ are defined. By $\overline{k_0}\subseteq k$ we denote the algebraic closure of $k_0$. Choice of a transcendental basis for $k_0/\Q$ provides an embedding ${k_0}\hookrightarrow \Ko$. Then we can extend it to an embedding $\overline{k_0}\hookrightarrow \Ko$ giving the base change $\iota: S_\Ko\to S_{\overline{k_0}}$ between Noetherian schemes. By \cref{lem: iso_fields}, we have $$\iota^{-1}\left(\Wlc_{(X, B)_{\overline{k_0}}}(\X/S, \B)_{\overline{k_0}}\right) = \Wlc_{(X, B)_\Ko}(\X/S, \B)_\Ko,$$ and the latter set is constructible by \cite[Theorem 4.3]{Ser26b}. Then $\Wlc_{(X, B)_{\overline{k_0}}}(\X/S, \B)_{\overline{k_0}}$ is constructible due to \cref{lem: constr_fields}. Therefore, the subset $\Wlc_{(X, B)}(\X/S, \B)\subseteq S$ is constructible again by \cref{lem: iso_fields} and \cref{lem: constr_fields}. The proof of constructibility for the subsets $\Iso_{(X, B)}(\X/S, \B)$, $\Iso^1_{(X, B)}(\X/S, \B)\subseteq S$ is analogues.
	\end{proof}
	
	Let $\Cc$ be a class or set of pairs. By $\|\Cc\|$ we denote the set of log isomorphism classes in $\Cc$.
	
	\begin{proposition}
		\label{prop: constr_and_finiteness}
		Let $(\X/S, \B)$ be a projective elementary family, and $(X, B)$ a klt wlc model. Let $W\subseteq S$ be a constructible subset such that $W\subseteq \Wlc_{(X, B)}(\X/S, \B)$. Suppose  for every $s\in W(k)$ the set $\Iso_{(\X_s, \B_s)}(\X/S, \B)\cap W$ is constructible. Then the set $\left\|\{(\X_s, \B_s): \ s\in W(k)\}\right\|$ is finite.
	\end{proposition}
	\begin{proof}
		Let us reduce the proof to the case when $k$ is uncountable. Suppose the proposition holds for all uncountable algebraically closed fields $K$. Let $K$ be an algebraic closure of $k((t))$, and $\iota: S_{K}\to S$ be the morphism induced by the field extension $K/k$. From \cref{lem: iso_fields} it follows that after the base change $\iota: S_{K}\to S$ we obtain a klt wlc model $(X, B)_K$, an elementary projective family $(\X/S, \B)_K$ and a constructible subset $W_K: = \iota^{-1}(W) \subseteq S_K$ satisfying the assumptions of the proposition. In addition, we have an inclusion:
		$
			\left\{\left(\X_s, \B_s\right)_K: \ s\in W(k)\right\}\subseteq 	\left\{\left((\X_K)_s, (\B_K)_s\right): \ s\in W_K(K)\right\}.
		$
		 Note that $k((t))$, as well as $K$, is uncountable by Cantor's diagonal argument. Hence, the following set is finite by our hypothesis:  $\left\|\left\{\left((\X_K)_s, (\B_K)_s\right): \ s\in W_K(K)\right\}\right\|.$ Again by \cref{lem: iso_fields}, we obtain the following isomorphism of sets which completes the reduction step: $$\left\|\left\{\left(\X_s, \B_s\right): \ s\in W(k)\right\} \right\| \cong \left\|\left\{\left(\X_s, \B_s\right)_K: \ s\in W(k)\right\}\right\|.$$
		
		Now, we can assume that $k$ is uncountable. Let $(X_{\al}, B_\al)$ be distinct pairs representing all classes in $\left\|\left\{\left(\X_s, \B_s\right): \ s\in W(k)\right\} \right\|$.  We claim that the latter set is at most countable. Indeed, let $g_s: (\X_s , \B_s)_\text{trm}\to (\X_s, \B_s)$ be a $\Q$-factorial terminalization \cite{BCHM10}. Choose an ample divisor $A_s$ on each pair $(\X_s, \B_s)$, put $D_s = g_s^* A_s$, and fix $o\in W(k)$. By \cref{def: 0-class}, there is an 0-equivalence $\varphi_s: (\X_s, \B_s)_\text{trm} \dashrightarrow (\X_o, \B_o)_\text{trm}$ which is isomorphism in codimension 1 by \cite{Kaw08}. Therefore, the cardinality of $\left\|\left\{\left(\X_s, \B_s\right): \ s\in W(k)\right\} \right\|$ is less or equal than the cardinality of $$\left\{\left[(\varphi_s)_*D_s\right]\in \text{NS}(\X_o): s\in W(k)\right\}.$$ The latter is at most countable, and the claim follows.  Let $W = \bigsqcup_i W_i$ be a finite decomposition into locally closed subschemes of $S$. Then we can write $W_i = \bigsqcup_\al \Iso_{(X_\al, B_\al)}(\X/S, \B)\cap W_i$ for distinct log pairs $(X_\al, B_\al)$ representing all sets in $\left\|\left\{\left(\X_s, \B_s\right): \ s\in W(k)\right\} \right\|$. Since $k$ is uncountable and each set $ \Iso_{(X_\al, B_\al)}(\X/S, \B)\cap W_i$ is constructible, one of $\Iso_{(X_\al, B_\al)}(\X/S, \B)$ must include a dense open subset of $W_i$. Then Noetherian induction on closed subsets of each $W_i$ finishes the proof.
	\end{proof}
		It turns out that the assumptions on singularities in the proposition above could not be weakened to log canonical, see the following example.
	\begin{example}
		Suppose $k = \Ko$. Let $C\subset \Pp^2_\Ko$ be a smooth cubic. Then $(X, B): = (\Pp^2_\Ko, C)$ is a projective plt 0-pair, which is not klt. Let $\Delta : = \{(s, s)\in \Pp^2_\Ko\times_\Ko C: \  s\in C\subset \Pp^2_\Ko\}$ be the diagonal, $\X: = \text{Bl}_\Delta(\Pp^2_\Ko\times_\Ko C)$ be its blowup, and $\B$ be the strict transform of $C\times_\Ko C$. Then $(\X/S, \B)$ is a projective elementary family over $S: = C$. By construction, $(\X_s, \B_s) = (\text{Bl}_s\Pp^2_\Ko, \widetilde{C}_s)$, where $\widetilde{C}_s$ is the strict transform of $C\subset \Pp^2_\Ko$ to the blow-up of $s\in C$. Furthermore, we have $S = \text{Wlc}_{(X, B)}(\X/S, \B)$, and $(\X_s, \B_s)\approx (\X_t, \B_t)$ if and only if there exists an automorphism $\varphi\in \text{Aut}(X, B)$ such that $\varphi(s) = t$. We claim that the group $\text{Aut}(X, B)$ is finite. If this is the case, then for every $s\in S(k)$ the set $\text{Iso}_{(\X_s, \B_s)}(\X/S, \B) = \text{Aut}(X, B)\cdot s$ is finite (constructible) but $\|(\X_s, \B_s): s\in S(k)\|$ is an uncountable set. 
		To prove the claim, we fix a flex point $o\in C$, and put $A = \Os_{C}(3o)$. Then  the restriction map  $\text{Aut}(\Pp^2_\Ko, C) \to \text{Aut}(C, \emptyset; A)$ is an isomorphism. Note that there is an isomorphism $\text{Aut}(C, \emptyset; A) \cong \text{Aut}(C, o)\ltimes E[3]$ where the latter group is finite (see \cite{Har77}).
	\end{example}
	
	\begin{lemma}[Shokurov's log geography]
		\label{lemma: log_geography}
		Let $(X, B)$ be a projective good wlc model, $Y$ be its log canonical model. Then there exists a decomposition $B = \sum_i r_i \Delta_i$ satisfying the following.
		\begin{itemize}
			\item 
				Real numbers $r_i\in [0, 1]$ are linearly independent over $\Q$, and $\sum_i r_i = 1$.
			\item 
				Each $\Delta_i$ is an effective $\Q$-divisor, and $\supp(\Delta_i) = \supp(B)$.
			\item 
				Each $(X, \Delta_i)$ is klt if and only if $(X, B)$ is klt.
			\item 
			Each $(X, \Delta_i)$ is a projective good wlc model whose log canonical model is $Y$.
		\end{itemize}
	\end{lemma}
	\begin{proof}
		In accordance with the convex approximation of $\R$-divisors II \cite[Proposition 11.47]{Kol26}, there exist (finitely many) real numbers $r_i\in [0, 1]$ linearly independent over $\Q$ and effective $\Q$-divisors $\Delta_i$ such that $B = \sum r_i \Delta_i$, $\text{supp}(\Delta_i) = \text{supp}(B)$ for all~$i$, and $\sum r_i = 1$. Moreover, each $(X, \Delta_i)$ is klt if and only if $(X, B)$ is klt. Let $g: X\to Y$ be the structure morphism, so we can write $K_X + B = g^*A + P$, where $A$ is an ample $\R$-Cartier $\R$-divisor, and $P$ is an $\R$-linear combination of principal divisors. Again by \cite[Proposition 11.46]{Kol26}, we can write $f^*A + P = \sum r_i f^*A_i$ for some ample $\Q$-Cartier $\Q$-divisors $A_i$. Therefore, $K_X + \Delta_i\sim_\Q f^*A_i$ for all $i$. In other words, each pair $(X, \Delta_i)$ is a good wlc model, and $g: X\to Y$ is its log canonical model.
	\end{proof}
	\begin{remark}
		\label{rem: log_geography}
		The first two statements in the proposition hold for any pair $(X, B)$ due to the convex approximation of $\R$-divisors I \cite[Proposition 11.46]{Kol26}, and the third one holds for any log pair.
	\end{remark}

	Let $(\X/S, \B)$ be an elementary family. For a given real number $\epsilon > 0$ define the following sets:
	 \begin{equation*}
	 	\begin{gathered}
	 		\epsilon\text{-}\textnormal{klt}(\X/S, \B) = \{s\in S: (\X_s, \B_s) \text{ is an $\epsilon$-klt log pair}\},\\
	 		\textnormal{klt}(\X/S, \B) = \{s\in S: (\X_s, \B_s) \text{ is a klt log pair}\},\\
	 		\wlc(\X/S, \B) = \{s\in S: (\X_s, \B_s) \text{ is a wlc model}\},\\
	 		_+\wlc(\X/S, \B) = \{s\in S: (\X_s, \B_s) \text{ is a good wlc model}\}.
	 	\end{gathered}
	 \end{equation*}
	\begin{proposition}
		\label{lemma: klt_0_families}
		Let $(\X/S, \B)$ be a projective elementary family, and $\B$ is an effective $\R$-divisor. Then for any real number $0< \epsilon< 1$ the following hold.
		\begin{enumerate}
			\item 
			If the set ${\textnormal{klt}_\epsilon(\X/S, \B)}$ is dense in $S$ then there exists a dense open subset $U\subseteq S$ such that $$(\X_U, \B_U) = (\X, \B)\times_S U\text{ is a klt log pair}.$$
			
			\item 
			If $(\X, \B)$ is a klt log pair, then the subset $\textnormal{klt}(\X/S, \B)\subseteq S$ is constructible.
			
			\item 
			The subset $\textnormal{klt}_\epsilon(\X/S, \B)\subseteq S$ is constructible.
		\end{enumerate}
	\end{proposition}
	\begin{proof}
		Let $\B = \sum r_i\Delta_i$ be a decomposition as in  \cref{lemma: log_geography}. Then ${\textnormal{klt}_\epsilon(\X/S, \Delta_i)} = {\textnormal{klt}_\epsilon(\X/S, \B)}$ for every index $i$.
		According to \cite[Proposition 2.4]{HX15}, for every $(\X, \Delta_i)$ there is a dense open $U_i\subseteq S$ such that $(\X, \Delta_i)_{U_i}$ is a klt log pair. Hence, $(\X, \B)_U$ is a klt log pair for $U = \bigcap_i U_i$. This proves (1).
		
		Next, we prove (2).  By Noetherian induction on closed subsets of $S$ in which $\textnormal{klt}(\X/S, \B)$ is dense, we can assume that $S = \overline{\textnormal{klt}(\X/S, \B)}$. Let $g: (\Y, \B_\Y) \to (\X, \B)$ be a log resolution, $K_\Y + \B_\Y = g^*(K_\X + \B)$. Shrinking the base, we may suppose that $(\Y/S, \B_\Y)$ is a projective elementary family, and $(\Y, \B_\Y)$ is log smooth over $S$. Then $g_s: (\Y_s, (\B_\Y)_s) \to (\X_s, \B_s)$ is a log resolution for all $s\in S$. By construction, $$g^*(K_\X + \B)|_{\Y_s}=(K_\Y + \B_\Y)|_{\Y_s} = (K_{\Y_s} +( \B_\Y)_s) = g_s^*(K_{\X_s} + \B_s).$$ This implies that $\text{dis}(\X_s, \B_s) = \text{dis}(\X, \B)$ for all $s\in S$, i.e. $(\X_s, \B_s)$ is a klt log pair for all $s\in S$.
		
		Finally, we prove (3). By Noetherian induction on closed subsets of $S$ in which $\textnormal{klt}_\epsilon(\X/S, \B)$ is dense, we can assume that $S = \overline{\textnormal{klt}_\epsilon(\X/S, \B)}$. Hence, (3) follows from (1) and (2).
	\end{proof}
	
	\begin{proposition}
		\label{prop: +wlc_constructible}
		Let $(\X/S, \B)$ be a projective elementary family, and $\B$ is an effective $\R$-divisor such that $(\X, \B)$ is a lc log pair. Then the subset $_+\wlc(\X/S, \B)\subseteq S$ is constructible.
	\end{proposition}
	\begin{proof}
	 As before, by Noetherian induction on closed subsets of $S$ in which $_+\wlc(\X/S, \B)$ is dense, we can assume that $S = \overline{_+\wlc(\X/S, \B)}$. Let $\B = \sum r_i\Delta_i$ be a decomposition as in  \cref{lemma: log_geography}. 
	 Choose a closed point $o\in _+\text{wlc}(\X/S, \B)$. From \cref{lemma: log_geography} it follows that $(\X_o, (\Delta_i)_o)$ is a good wlc model for every index $i$. This implies the following inclusion:
		\begin{equation*}
			_+\text{wlc}(\X/S, \Delta_i) \supseteq _+\text{wlc}(\X/S, \B) \ \text{ for all } i.
		\end{equation*}
	Let $g: (\Y ,\B_\Y)\to (\X, \B)$ be a log resolution, and $f: \X\to S$ the structure morphism.  Shrinking the base, we may suppose that $(\Y/S, \B_\Y)$ is a projective elementary family, and $(\Y, \B_\Y)$ is log smooth over $S$. 
	Then we write $\B_\Y = \sum r_i (\B_\Y)_i$ where $K_\Y + (\B_\Y)_i = g^*(K_\X + \Delta_i)$. Next, for each $\Q$-divisor $(\B_\Y)_i$ we can write $(\B_\Y)_i = \Delta_i - E_i$ where  each $E_i$ is supported on the exceptional divisors, and $\Delta_i, E_i\ge 0$. By \cite[Corollary 1.4]{HMX18}, for every closed point $s\in S$ the canonical morphism $$(f\circ g)_*\Os_\Y\left(m\left(K_{\Y} + \Delta_i\right)\right)\otimes \kappa(s)\to H^0\left(\Y_s, \Os_{\Y_s}\left(m\left(K_{\Y_s} +( \Delta_i)_s\right)\right)\right)$$ 	is surjective for every $m\in \mathbb{Z}_{>0}$ such that $m\Delta_i$ is integral. From now and on, fix $m\in \mathbb{Z}_{>0}$ such that $m\Delta_i$ is integral for all $i$. Put $L_i = \Os_\X\left(m\left(K_{\X} + \Delta_i\right)\right)$. Recall that $K_\Y + \Delta_i = g^*(K_\X + \Delta_i) + E_i$. Hence, the canonical morphism $(f_*L_i)\otimes \kappa(s)\to H^0\left(\X_s, L_i|_{\X_s}\right)$ is surjective as well. Set $Z_i = \text{supp coker}\left(f^*f_*L_i\to L_i\right)$ be the relative base locus. Then every $Z_i$ is vertical over $S$. Indeed, for any $s\in \, _+\text{wlc}(\X/S, \Delta_i)$ there is a section $\sigma\in H^0(\X_s, L_i|_{\X_s})$ missing $Z_i\cap \X_{s}$. As we proved, the section $\sigma$ can be lifted to $L_i$. Therefore, every $L_i$ is semiample over $S\smallsetminus \bigcup_i f(Z_i)$. This concludes the proof.
	\end{proof}
	
	Although, the following proposition will not be used later in the paper, it could be of independent interest and useful in other contexts.
	\begin{proposition}
		Let $(\X/S, \B)$ be a projective elementary family, $\B$ is an effective $\R$-divisor, and $(\X, \B)$ is a lc log pair. Suppose the conjecture on existence of minimal models holds for $(\X, \B)/S$. Then the subset $\wlc(\X/S, \B)\subseteq S$ is constructible.
	\end{proposition}
	\begin{proof}
		Let $\varphi: (\X, \B)/S \dashrightarrow (\X', \B')/S$ be a minimal model for $(\X, \B)/S$, that is, $(\X',\B')$ is a lc log pair, projective over $S$, and $K_{\X'} + \B'$ is nef over $S$. Also, the definition of a minimal model assumes that there is a common resolution $W$  of $(\X, \B)$ and $(\X', \B')$ such that $p^*(K_\X + \B) = q^*(K_{\X'} + \B') + E$ for some effective $q$-exceptional $\R$-divisor $E$, here $p: W\to \X$ and $q: W\to \X'$ are projective birational contractions. Then, the set $Z = \bigcup_{E_i\in \text{supp}(E)}f(E_i)$ is closed because $W/S$ is a projective morphism. Note we can assume that $(W/S, E)$ is a projective elementary family and $p_s^*(K_{\X_s} + \B_s) = q^*(K_{\X'_s} + \B'_s) +  E_s$ for all $s\in U$ for some dense open $U\subseteq S$. By the same argument as in the proof of \cref{lemma: klt_0_families}, $(\X_s, \B_s)$ is lc for every $s\in S$ after shrinking the base. Hence, $\wlc(\X/S, \B) = \{s\in S: K_{\X_s} + \B_s \text{ is nef }\} = S\smallsetminus Z$ is open. By Noetherian induction, we conclude that $\wlc(\X/S, \B)$ is constructible.
	\end{proof}

	\section{Main results}
	\begin{proposition}
		\label{lem: rigidity_small}
		Let $(\X, \B)/S$ be a projective klt wlc model over a smooth variety $S$. Suppose there is a klt wlc model $(X', B')\times_k S/S$ in the 0-class of $(\X, \B)/S$. Then $(\X, \B)/S \approx (X, B)\times_k S/S$ for some klt wlc model $(X, B)$. 
	\end{proposition}
	\begin{proof}
		After replacing $(X', B')$ by a general fiber of $(\X, \B)/S$, we can assume that $(X', B')$ is projective. Let $({X'}_\text{trm}, {B'}_\text{trm})$ and $(\X_\text{trm}, \B_\text{trm})$ be $\Q$-factorial terminalizations \cite{BCHM10} of $({X'}, {B'})$ and $(\X, \B)$ respectively. Let us remark that $(X', B')\times_k S$ has (at worst) terminal singularities because $S$ is smooth. According to \cite{Kaw08} and \cref{lemma: log_geography}, a birational map $\varphi: ({X'}_\text{trm}, {B'}_\text{trm})\times_k S \dashrightarrow (\X_\text{trm}, \B_\text{trm})$ is small. Choose a general ample  Cartier divisor $\A$ on $\X_\text{trm}/S$ such that $(\X_\text{trm}, \B_\text{trm} + \epsilon \A)$ and $(X_\text{trm}\times_k S, \ B_\text{trm}\times_k S + \varphi_*^{-1}\A)$ are klt log pairs for some real number $0<\epsilon\ll 1$.  Then, $\varphi$ decomposes into a finite sequence of flips over $S$ with respect to $(K_{X'_\text{trm}} + B'_\text{trm})\times_k S + \epsilon\varphi_*^{-1}\A$ again by \cite{Kaw08}. From  the rigidity lemma \cite[Lemma 1.15]{Deb01} it follows that each flip preserves the triviality over $S$ (see \cite[Lemma 5.2(2)]{FS25}). In other words, there is a trm log pair $(X_\text{trm}, B_\text{trm})$ with an isomorphism $(\X_\text{trm}, \B_\text{trm})/S\approx (X_\text{trm}, B_\text{trm})\times_k S/S$. Then we conclude that $(\X, \B)/S \approx (X, B)\times_k S/S$ for some klt wlc model $(X, B)$  again by the rigidity lemma (see \cite[Lemma 5.2(1)]{FS25}).
	\end{proof}
	\begin{remark}
		There is a more general rigidity result \cite[Theorem 30]{Kol26} about when a rational product structure can be uniquely extended to a regular product structure.
	\end{remark}

	\noindent For an elementary family $(\X/S, \B)$ and a proper log pair $(X, B)$ we set
	$$_+\Wlc_{(X, B)}(\X/S, \B):= \, _+\wlc(\X/S, \B)\cap \Wlc_{(X, B)}(\X/S, \B).$$ The sets $_+\Iso_{(X, B)}(\X/S, \B)$ and $_+\Iso^1_{(X, B)}(\X/S, \B)$ are defined similarly. Let us recall that a \textit{cover} is a finite surjective morphism; an \textit{\'etale base change} $\iota: T\to S$ is a morphism between smooth varieties such that $\iota(T)\subseteq S$ is dense open and $\iota: T\to \iota(T)$ is an \'etale cover. A cover is called \textit{quasi-\'etale} if it is \'etale in codimension 1.
	\begin{theorem}
		\label{th: main_constr}
		Let $(\X/S, \B)$ be a projective elementary family, and $(X, B)$ a wlc klt model. Then the following subsets are constructible: $$_+\Iso_{(X, B)}(\X/S, \B), \, _+\Iso_{(X, B)}^1(\X/S, \B), \,  _+\textnormal{Wlc}_{(X, B)}(\X/S, \B)\subseteq S.$$
	\end{theorem}
	\begin{proof}
		 Put $_+\Wlc_{(X, B)}(\X/S, \B) := W$. Let $\mathcal{S}$ be the set of all closed subsets $S^\prime \subseteq S$ such that $S^\prime \cap W$ is dense in $S^\prime$. Note that $\mathcal{S}$ is a well-founded set with respect to inclusion as $S$ is Noetherian. Then the set of all minimal elements for $\mathcal{S}$ is $\{\{s\}: s\in W(k)\}$. We say that a property $P(S^\prime)$ holds for $S^\prime \in \mathcal{S}$ if the set $W\cap S^\prime$ is constructible. We shall proceed by Noetherian induction on $\mathcal{S}$. Clearly, $P(\{s\})$ holds for all $s\in W(k)$ (Noetherian induction base). Let $S^\prime\in \mathcal{S}$ be a non-minimal element of $\mathcal{S}$. Assume that $P(S^{\prime\prime})$ is true for all subsets $S^{\prime\prime}\subset S^\prime$ such that $S^{\prime\prime}\in \mathcal{S}$ (Noetherian induction step). We claim that $P(S^\prime)$ holds. To prove the claim, it suffices to show that there is a generically finite base change $U\to S^\prime$ such that $(\X_U, \B_U) = (\X, \B)\times_S U \approx (X', B')\times_\Ko U$ over $U$ for some klt wlc model $(X', B')\sim_\text{bir} (X, B)$. Then the statement follows from Noetherian induction. It is readily seen that the problem can be reduced to the case when $(\X_{S^\prime}/S^\prime, \B_{S^\prime})$ is a projective elementary family of good klt wlc models according to \cref{lemma: klt_0_families}. To simplify the notation, we write $(\X/S, \B)$ for the corresponding new family. To be definite, we fix a point $o\in W(k)$. After each generically finite base change $\iota: T\to S$, we replace $o$ by a point in $\iota^{-1}(W)(k)$. From now and on, we denote $(\X_o, \B_o)$ by $(X, B)$.
		 
		 \textbf{Step 1.} We reduce the proof to the case when $\B$ is a $\Q$-divisor by applying Shokurov's log geography. Suppose the theorem holds for all projective elementary families with $\Q$-coefficients. Let $\B = \sum r_i \Delta_i$ be a decomposition as in \cref{lemma: log_geography}.  From \cref{lemma: log_geography} it follows that $(\X_s,( \Delta_i)_s)$ is a good klt wlc model for all $S(k)$. Set $\Delta^0_i:= (\Delta_i)_o$. Therefore, we obtain $$ W\subseteq \ \text{Wlc}_{(X, \Delta^0_i)}(\X/S, \Delta_i) \ \text{ for all indices } i.$$ Hence, after shrinking the base, we can assume that $S = \ \text{Wlc}_{(X, \Delta^0_i)}(\X/S, \Delta_i)$ by the hypothesis. From now and on, we fix one index $i$. Next, according to \cref{prop: constr_and_finiteness} there is a finite decomposition $S = \bigsqcup_\al\text{Iso}_{(X_\al, B_\al)}(\X/S, \Delta_i)$ corresponding distinct klt wlc models $(X_\al, B_\al)$. Again by our hypothesis, every locus $\text{Iso}_{(X_\al, B_\al)}(\X/S, \Delta_i)$ is constructible.  Consequently, after an \'etale base change, there is an isomorphism $\Phi: (X, B_\text{red})\times_k S \to (\X/S, \B_\text{red})$ over $S$, where $(X, B) = (X_\al, B_\al)$ for some index $\al$. Then $\Phi_o = \text{id}_X \in \text{Aut}(X, B_\text{red})$ according to our convention. By construction, the log isomorphism $\varphi:= (\Phi_o\times_k \text{id}_S)\circ \Phi^{-1}: (\X, \B_\text{red}) \to (X, B_\text{red})\times_k S$ over $S$ provides $(\varphi_s)_*(\B_s) = B'$ for every $s\in S$, as  required.
		 
		 \textbf{Step 2.} We show that after an \'etale base change there is a projective contraction $g: \X\to Y\times_k S$ over $S$ where $Y$ is the log canonical model for $(X, B)$, when $\B$ is a $\Q$-divisor.
		 Let $L = \Os_{\X}\left(m(K_{\X} + \B)\right)$ for $m\gg 1$. From \cite[Corollary 1.4]{HMX18} it follows that, after shrinking the base, there exists a projective contraction $g: \X/S\to \Y/S$ such that $L\approx_S g^*(\A)$ for some ample invertible sheaf $\A$ on $\Y/S$,  cf. the proof for \cref{prop: +wlc_constructible}. Put $A = \A_{g_{o}}$ and $Y = \Y_o$. Then the set $\Iso_{(Y, \emptyset; A)}(\Y/S, \emptyset; \A)$, which corresponds to isomorphic polarized pairs, is constructible by \cite[Theorem 4.12]{Ser26b}. Since the log canonical model is unique up to log isomorphism for 0-equivalent pairs, we have $W\subseteq \Iso_{(Y, \emptyset; A)}(\Y/S, \emptyset; \A)$. Thus, the latter is dense in $S$. Hence, after an \'etale base change we can assume that $(\Y, \emptyset; \A)/S = (Y, \emptyset; A)\times_k S /S$.
		  
		  \textbf{Step 3.} We prove that there is an \'etale base change with an automorphism of $Y\times_k S/S$ such that for every $s\in W(k)$ there is a 0-equivalence $\varphi_s: (\X_s, \B_s) \dashrightarrow (X, B)$ over $Y$. By construction, for every point $s\in W(k)$ there is an automorphism $\theta_s: (Y, \emptyset; A)\times_k s\to (Y, \emptyset; A)\times_k o$ and a 0-equivalence $\varphi_s: (\X_s, \B_s) \dashrightarrow (X, B)$ such that $g_o\circ \varphi_s = \theta_s\circ g_s $.
		  Let $H$ be the group scheme representing the functor $T \mapsto \text{Aut}(Y_T, \emptyset; A_T)$. It is well-known that $H$ is of finite type \cite[Theorem 2.16]{Bri18}. Let $H_i$ be all irreducible components of $H$. By $AH_i$ we denote the closure of $\{(\theta_s, s) : s\in W(k)\}$ in $H_i\times_k S$. According to Chevalley's theorem $AH:=\bigcup_i \text{pr}\left(AH_i\right)\subseteq S$ is a constructible subset. Note $AH\subseteq S$ is dense because $W(k)\subseteq AH$. Hence, one of $\text{pr}(AH_i)$ must include a dense open subset of $S$. After an \'etale base change, we can assume that there exists a section for $S\to AH_i$ that induces an automorphism $\Theta: (Y, \emptyset; A)\times_k S \to (Y, \emptyset; A)\times_k S$ over~$S$. Therefore, after replacing $g$ by $\Theta\circ g$ we can assume that for all $s\in W(k)$ given 0-equivalence $\varphi_s: (\X_s, \B_s) \dashrightarrow (X, B)$ is defined over~$Y$.
 		 
		 \textbf{Step 4.} We prove that after an \'etale base change the induced family $(\X/S, \B)_{\overline{\eta}}$ over the geometric generic point $\overline{\eta}$ of $Y$ is trivial, when $\B$ is a $\Q$-divisor.
		 Let $\text{pr}: Y\times_k S\to Y$ be the projection. Note that $g: (\X, \B)\to Y\times_k S$ is a projective lc-trivial fibration by construction. Let $B_{Y\times_k S}$ denote its boundary part. After shrinking the base, we can assume that every prime component of $B_{Y\times_k S}$ is horizontal over $S$. Then for the generic point $\eta\in Y$ the induced morphism $g_{K}: (\X, \B)_{K}\to S_{{K}}$ is also a projective lc-trivial fibration, where $K$ is an algebraic closure of $ \kappa(\eta)$.  Denote by $\iota: S_{K}\to S$ the induced morphism. Set $W':= \Wlc_{(X, B)_{K}}(\X/S, \B)_{K}$. Recall that for every point $s\in W(k)$ there is a 0-equivalence $\varphi_s: (\X_s, \B_s) \dashrightarrow (X, B)$ over $Y$ by the previous step. Consequently, $\iota^{-1}(W)\subseteq W$. Let us note that $\iota^{-1}(W)$ is Zariski dense in $S_{K}$ because $W(k)$ is Zariski dense in $S$. Furthermore, there exist finitely many pairs $(F_\al, \Delta_\al)$ over $K$ such that $W' = \bigsqcup_\al \text{Iso}_{(F_\al, \Delta_\al)}(\X/S, \B)_{K}$ according to \cref{prop: constr_and_finiteness} and \cref{th: 0-pairs_constr}. Hence, for some index $\al$ the subset $\Iso_{(F_\al, \Delta_\al)}(\X/S, \B)_K$ includes a dense open subset $U\subseteq S_{K}$. To ease the notation, denote  this pair by $(F, \Delta)$. Note the boundary part $B_{S_K}$ for $(\X/S, \B)_{K}$ is zero because no components of $B_{Y\times_k S}$ dominates $Y$. In addition, the moduli part $M_{S_K}\sim_\Q 0$ because $(\X/S, \B)_K$ is isotrivial over $U$. Therefore, from \cite[Theorem 4.7]{Amb05} it follows that there is a quasi-\'etale cover  $T_K\to S_{K}$ from a regular variety such that that $(\X, \B)_{T_K} = (F, \Delta)\times_{K} T_K$. By the Zariski--Nagata theorem on purity of branch locus \cite[\href{https://stacks.math.columbia.edu/tag/0BMB}{Tag 0BMB}]{Stacks}, the cover $T_{K}\to S_{K}$ is \'etale because $S_{K}$ is regular. Hence, the cover $T_{K}\to S_{K}$ is induced by some \'etale cover $T\to S$  as the natural functor $\text{F\'et}_{S}\to \text{F\'et}_{S_K}$ is an equivalence, see \cite[Theorem 1.1]{Lan24} and references therein. From now and on, we assume $(\X/S, \B)_{K}$ is trivial.
		 
		 \textbf{Step 5.} We show that after an \'etale base change of $S$, the family $(\X/S, \B)_\eta$ becomes trivial. To prove this, we follow the descent argument in \cite[Theorem 4.7 (4)]{Amb05}. According to the previous step, the family $(\X/S, \B)_K$ is trivial. Hence we can find a finite Galois extension $L/\kappa(\eta)$ such that $(\X/S, \B)_L$ is trivial, i.e. there exists an isomorphism $\Phi: (\X/S, \B)_L \to (F, \Delta)_L\times_L S_L$ over $S_L$. Recall that $(X, B) = (\X_o, \B_o)$. Then we can assume that $(F, \Delta) = (X, B)_L$ and $\Phi_{o_L} = \text{id}_F$ for some $o_L\in S_L(L)$ that maps to $o\in S(k)$.  Next, a Galois automorphism $\sigma\in G:= \text{Gal}(L/\kappa(\eta))$ induces natural automorphism $\beta_\sigma := \text{id}_F\times_L \text{Spec}(\sigma^{-1})$ which may not be equal to $\sigma': =\Phi\circ\sigma \circ \Phi^{-1}$. Hence, we obtain \textit{a priori} non-trivial automorphisms $\sigma'\circ \beta_\sigma^{-1}$ of $(F, \Delta)\times_L S_L$ over $S_L$. In other words, there is a morphism $t_\sigma: S_L\to \text{Aut}(F, \Delta)$. As $\Phi_{o_L} = \text{id}_F$, each $t_\sigma$ factors through $A_L: = \text{Aut}^0(F, \Delta)$ that is an abelian variety \cite[Proposition 4.6]{Amb05}.  It is clear that the set of morphisms $A(S_L) := \text{Mor}_{L}(S_L, A_L)$ forms a group. Note that there is a natural left action of $G$ on  $A(S_L)$ given by $\sigma\cdot a(s) = \beta_\sigma(\_, s)\circ a(s)\circ \beta_\sigma^{-1}(\_, s)$ for every $s\in S_L(L)$.  Then the sections $t_\sigma$ satisfy the 1-cocycle relation:
		 \begin{equation*}
		 	t_{\sigma \tau}= t_\sigma\circ (\sigma \cdot t_\tau) \ \text{ for all } \sigma,\tau\in G.
		 \end{equation*}
		 Since $A$ is abelian, we shall use the additive notation. Set $t^0: = \sum_{\tau\in G} t_\tau$, $n := |G|$. Hence, for every $t_\sigma$ the 1-cocyle relation implies  that $nt_\sigma= t^0 - \sigma\cdot t^0$, that is, $\{\sigma\mapsto t_\sigma\}$ is a a torsion element in the Galois cohomology group $H^1(G, A(S_L))$ (cf. \cite[Lemme 7]{Ser58}). Hence, after an \'etale base change of $T_L\to S_L$, we can trivialize all $s_\sigma$.
		 Recall that $K$ is an algebraic closure for $\kappa(\eta)$, and $K\supseteq L$. Set $T_K := T_L\times_{S_L} S_K$. Again by \cite[Theorem 1.1]{Lan24}, the \'etale cover $T_{K}\to S_{K}$ is induced by some \'etale cover of $S$. Consequently, $(\X/S,\B)_{\kappa(\eta)}$ is trivial after an \'etale base change, as claimed.
		 
		 \textbf{Step 6.} Finally, we finish the proof by applying rigidity results for birationally trivial families. 
		 Then by the standard descent argument, we can find a dense open subset $U\subseteq Y$ with open dense subset $V\subseteq X$ such that $(\X, \B)|_{U\times_k S} = (V, B|_{V})\times_k S$.  Let us remember that $(X, B) = (\X_o, \B_o)$ for some $o\in W(k)$. Hence, $(\X,\B)/S$ is a klt wlc model which is 0-equivalent to $(X, B)\times_k S$.  By \cref{lem: rigidity_small}, the family $(\X/S, \B)$ is trivial. This finishes the proof.
		 
	\end{proof}

	In an earlier version of the paper, the following rigidity result was used in the proof of \cref{th: main_constr}. Despite of its replacement by a Galois cohomology argument, the proposition can be of independent interest and useful in other contexts.
	\begin{proposition}
		\label{lem: rigidity_quotient}
		Let $(\X/S, \B)$ and $(\X'/S, \B')$ be projective elementary families such that
		\begin{itemize}
			\item 
			$(\X', \B')/S\approx (X', B')\times_k S/S$ for some pair $(X', B')$,
			\item 
			$(\X/S, \B) \approx (\X'/S, \B')\big/ G$ for some finite group $G$.
		\end{itemize}
		Then $(\X/S	, \B)$ is \'etale locally trivial.
	\end{proposition}
	\begin{proof}
		Let $H$ be the group scheme representing a natural functor $T\mapsto \text{Aut}(X'_T, B'_T)$. Recall that every finite group in characteristic zero is reductive by Maschke's theorem. Then by \cite[Theorem 2]{Bri22}, the functor $\textbf{Hom}_\text{gr}(G, H)$ of morphisms of group schemes is represented by a smooth scheme $M$, and the morphism $H\times_k M\to M\times_k M$ given by $(h, \rho)\mapsto (h\rho h^{-1}, \rho)$ is smooth. Since $G$ acts on $\X'\approx X'\times_k S$ over $S$, there is a morphism $\Phi: S\to M$ sending every closed point $s$ to a group homomorphism $\rho_s: G \to \text{Aut}(X')$. Then the image of $\Phi$ is included in a connected component of $M$ because $S$ is irreducible. By \cite[Remark 6.4]{Bri22}, the connected components of $M$ are precisely $H^0$-orbits of $k$-points, where $H^0$ is the connected component of the identity in $H$. Hence, for every closed point  $s\in S$ there is $h_s\in H^0$ such that $\rho_s = h_s\rho_{o}h_s^{-1}$ for some fixed $o\in S(k)$. This implies that  $(\X/, \B)$ is \'etale locally trivial.
	\end{proof}
	\begin{remark}
		More generally, the same argument works for linearly reductive algebraic groups $G$ and smooth groups $H$ in arbitrary characteristic via \cite[Theorem 6.3]{Bri22}.
	\end{remark}

	\section{Applications}
	\label{sec: applications}
	\begin{theorem}
		\label{th: boundedness_and_finiteness}
		Let $\Cc$ be a collection of 0-equivalent projective good klt wlc models. Then the set $\|\Cc\|$ of log isomorphism classes in $\Cc$  is finite if and only if $\Cc$ admits a bounded polarization.
	\end{theorem}
	\begin{proof}
			According to \cref{prop: bounded_subclass}, there exists finitely many projective elementary families $(\X^{(j)}/S^{(j)}, \B^{(j)})$ such that for every pair $(X_\alpha, B_\alpha)$ in $\Cc$ there is an index $j$ together with a closed point $s_\alpha\in S^{(j)}$, and an isomorphism $\varphi_\al: (X_\al, B_\al) \to (\X_{s_\al}^{(j)}, \B_{s_\al}^{(j)})$. Put $W_j = _+\Wlc_{(X, B)}(\X^{(j)}/S^{(j)}, \B^{(j)})$. In accordance with our main result, \cref{th: main_constr}, for every index $j$  the set $_+\Iso_{(\X_s, \B_s)}(\X^{(j)}/S^{(j)}, \B^{(j)})$ is constructible for every $s\in W_j(k)$, and $W_j$ is so. Then the theorem follows directly from \cref{prop: constr_and_finiteness}.
	\end{proof}

	Let us recall that for a proper surjective morphism $X/Y$ between varieties by $\text{NS}(X/Y)$ we denote the group of Cartier divisors on $X$ modulo the numerical equivalence over $Y$, and $X_\eta$ be the generic fiber. Then the group $\text{Aut}(X/Y)$ has a natural right action on $\text{NS}(X/Y)$ defined by pullbacks. Moreover, $\Aut(X/Y)$ also naturally acts on $\text{NS}(X_\eta)$ using the restriction map $\text{NS}(X/Y)  \to \text{NS}(X_\eta)$. By $[A_\eta]$ we denote the restriction of $[A]\in \text{NS}(X/Y)$ to $\text{NS}(X_\eta)$. The volume of a Weil divisor $D$ on a proper variety $X$ we define by the following formula:
	\begin{equation*}
		\text{vol}(D) = \underset{m\to +\infty}{\lim\sup} \frac{h^0(X, \Os_X(mD))}{m^d/d!},
	\end{equation*}
	where  $\Os_X(mD)$ is the divisorial sheaf associated with the Weil divisor $mD$, and $d:= \dim X$.
	\begin{theorem}
		\label{th: finiteness_big_nef}
		Let $(X, B)$ be a projective good klt wlc model, $Y$ its log canonical model, and $X_\eta$ the generic fiber of $X/Y$.  Then for every integer $v\in \mathbb{Z}$, the number of orbits for the natural action of $\textnormal{Aut}(X/Y, B)$ on the following set is finite:
		\begin{equation*}
			\{ \ [A_\eta]\in \textnormal{NS}(X_\eta): \textnormal{vol}(A_\eta) = v,\  A \textnormal{ is a Cartier divisor on $X$ which is big and nef over $Y$ } \}.
		\end{equation*}
	\end{theorem}
	\begin{proof}
		\textbf{Step 1.} We reduce the proof to the case when $B$ is a $\Q$-divisor. Suppose the theorem holds for all projective good klt wlc models with $\Q$-coefficients. Let $B = \sum r_i\Delta_i$ be a decomposition as in \cref{lemma: log_geography}. Then $\text{Aut}(X/Y, B) = \bigcap_i \text{Aut}(X/Y, \Delta_i)$. In fact, the index of $\text{Aut}(X/Y, B)$ in $\text{Aut}(X/Y, \Delta_i)$ is bounded by $\#\supp(B)!$. Hence, the finiteness of orbits for $\text{Aut}(X/Y, \Delta_i)$ implies the finiteness of orbits for $\text{Aut}(X/Y, B)$. From now and, we assume $B$ is an (effective) $\Q$-divisor. 
		
		\textbf{Step 2.}\label{Step 2} We prove finiteness of orbits for the natural action of $\text{Aut}(X/Y, B)$ on the set $$\{[A]\in \textnormal{NS}(X_\eta): \textnormal{vol}(A|_F) = v,\  A \textnormal{ is a Cartier divisor on $X$ which is ample over } Y\}.$$ Let $A_\al$ be (possibly non-effective) divisors representing all classes in the given set. By $\eta\in Y$ we denote the generic point. Then by \cite[Theorem 3.6]{Zhu26} and the strong Lefschetz principle, there is an integer $l\in \Zi_{>0}$ satisfying the following.
		\begin{itemize}
			\item 
				For every $A_\al$ there exists a  $\Q$-factorial projective log pair $(X_\al, B_\al)$ with a 0-equivalence $\psi_\al: (X, B)\dashrightarrow(X_\al, B_\al)$ which is an isomorphism in codimension 1. In addition, there is an effective divisor $J_\al$ such that $(J_\al)_\eta \sim l(A_\al')_\eta$ where $A'_\al := (\psi_\al)_* A_\al$.
			
			\item 
				The corresponding set $\Dd$ of all pairs $(X_\al, B_\al + J_\al)$ admits a bounded polarization.
		\end{itemize}
		By \cref{prop: bounded_subclass}, there are finitely many elementary projective families $(\X^{(j)}/S^{(j)}, \B^{(j)} + \mathcal{J}^{(j)})$ such that for every pair $(X_\alpha, B_\alpha + J_\al)$ in $\Dd$ there is an index $j$ together with a closed point $s_\alpha\in S^{(j)}$ and an isomorphism $\varphi_\al: (X_\al, B_\al + J_\al) \to (\X_{s_\al}^{(j)}, \B_{s_\al}^{(j)} + \mathcal{J}^{(j)}_{s_\al})$. Note that  $(\varphi_\al)_*{J}_\al = \mathcal{J}^{(j)}_{s_\al}$ for all $\al$ because $\text{coeff} (B_\al) \cap\text{coeff} (J_\al) =~\emptyset$. Applying Noetherian induction, we reduce the problem to one projective elementary family $(\X/S, \B + \mathcal{J})$ such that $W:= _+\Wlc_{(X, B)}(\X/S, \B)$ is dense in $S$.
		
		Let $\Cc$ denote the set  of underlying good wlc klt pairs $(X_\al, B_\al) \in \Dd$. Due to \cref{th: boundedness_and_finiteness}, the set $\|\Cc\|$ of log isomorphism classes in $\Cc$ is finite. Then from \cref{th: main_constr} it follows that after an \'etale base change we can assume $(\X, \B) = (X', B')\times_k S$ for some klt wlc model $(X', B')\in \Cc$. To simply the notation, we re-denote $\varphi_\al\circ\psi_\al$ by $\psi_\al$. Let $o\in S$ be a closed point such that $o = s_{\alpha_0}$ for some index $\al_0$. Therefore, $\mathcal{J}_s\sim_\text{alg} \mathcal{J}_o$ for all $s\in S(k)$. By construction, we have
		\begin{equation*}
			(lA'_\al)_\eta \sim (\mathcal{J}_{s_\al})_\eta \sim_\text{alg} (\mathcal{J}_{o})_\eta \sim (lA'_{\al_0})_\eta.
		\end{equation*}
		Note that pseudo-automorphisms preserve linear, algebraic, and numerical equivalences for Cartier divisors when their image is Cartier as well. Hence, we obtain that $\theta_\al := (\psi_{\al_0})^{-1}\circ\psi_\al \in \text{PsAut}(X/Y, B)$ and $(\theta_\al)_* (A_\al)_\eta \sim_\text{alg} A_{\al_0}$. Then $(\theta_\al)_\eta \in \text{Aut}(X_\eta, B_\eta)$ but it is possible that $\theta_\al \notin \text{Aut}(X/Y, B)$ as $(\theta_\al)_* (A_\al)$ is merely movable over~$Y$.
		
		 By construction, there are vertical over $Y$ Cartier divisors $V_\al$ such that $(\theta_\al)_* lA_\al \equiv_Y M_0 + V_\al$ where $M_0:= (\psi_{\al_0})^{-1}\mathcal{J}_o$. Let $E_i$ be prime vertical Cartier divisors generating the kernel of the restriction map $\text{NS}(X/Y)\to \text{NS}(X_\eta)$. Let $H_i\subset Y$ be a prime Cartier divisor such that $H_i\supseteq f(E_i)$, here $f: X\to Y$ is the structure map. Then for every $\al$ there is $n(\al)\in \Zi_{>0}$ such that $V'_\al = V_\al + \sum_i n(\al)f^*H_i$ is effective. Hence, $(\theta_\al)_*(lA_\al)\equiv_Y M_0 + V'_\al$ which is an effective Cartier divisor.
		 
		  Set $M_\al = M_0 + V'_\al$. Let us remember that for all $\al$ we have $\text{supp}(M_\al)\subseteq \bigcup_i\{E_i\}\cup \supp(M_0)$. Then for sufficiently small $0< \epsilon_\al\ll 1$, the maps $(\theta_\al)^{-1}: (X, B) \dashrightarrow (X, B)$ are ample models over $Y$ for $(X, B + \epsilon_\al M_\al)$. According to Shokurov's log geography \cite[Lemma 5.3]{JJZ26}, there are at most finitely many ample models for all $(X, B+\epsilon M_\al)$. Reducing the question to isomorphic ample models, we obtain that the maps $\theta_{\al\beta}:= (\theta_\al)^{-1}\circ\theta_\beta\in \text{PsAut}(X/Y, B)$ are log isomorphisms such that $(\theta_{\al\beta})_*A_\beta \equiv_Y A_\al + V_\al$. This finishes the proof of the finiteness of ample orbits.
		
		\textbf{Step 3.} We prove the theorem in general. Let $D_\al$ be (possibly non-effective) divisors representing all classes in the given set. Then by \cite[Theorem 2.2.4]{Fuj09} there exists a positive integer $l$ such that $\Os_X(lD_\al)$ is globally generated over $Y$ for all $\alpha$. Let $f_\al: X\to X_\al$ be the projective contraction over $Y$ given by $\Os_X(D_\al)$. Hence, we can write $lD_\al \equiv_Y \widetilde{D}_\al$ for some effective divisor $\widetilde{D}_\al \equiv_Y lD_\al$. Let  $A_\al$ be a relatively ample effective divisor such that $\widetilde{D}_\al = f_\al^*(A_\al)$. Put $B_\al = (f_\al)_* B$. Then the log pair $(X_\al, B_\al)$ is a projective klt wlc model which is 0-equivalent to $(X, B)$. By construction $\text{vol}(A_\al|_{F{\al}}) = l^rv$, where $F_\al$ is general fiber of $X_\al/Y$ and $r = \dim X - \dim Y$. Therefore, by \cite{Zhu26}, all pairs $(X_\al, B_\al)$ admits a bounded polarization. Hence, the set $\|\{(X_\al, B_\al)\}\|$ of log isomorphism classes is finite by \cref{th: boundedness_and_finiteness}. So, we can reduce the proof to the case when all targets $(X_\al, B_\al)$ are log isomorphic. Denote this target by $(X', B')$. Now, we can invoke the claim proven in \hyperref[Step 2]{Step 2} to reduce the proof to the case when all $[A_\al]\in \text{NS}(X'_\eta)$ belong to the same orbit of $\text{Aut}(X'/Y, B')$. Denote an effective generator of the corresponding orbit by $A'$. Therefore, for $0< \epsilon \ll 1$, all log pairs $(X, B + \epsilon \widetilde{D}_\al)$ are wlc klt models for $(X', B' + \epsilon A')$ over $Y$. \cite{BCHM10}, there are at most finitely many of them up to (natural) log isomorphism over $Y$ (cf. \cite[Lemma 2.8]{FHS25}).
	\end{proof}

	Similarly, 	for a proper surjective morphism $X/Y$ between varieties by $\text{NS}_W(X/Y)$ we denote the group of Weil divisors on $X$ modulo the numerical equivalence over $Y$, which is defined when $X$ is normal and $\Q$-factorial. Let $X_\eta$ be the generic fiber. Then the group $\text{PsAut}(X/Y)$  of isomorphisms in codimension 1 has a natural right action on $\text{NS}_W(X/Y)$ defined via inverse strict transforms.  Moreover, $\text{PsAut}(X/Y)$ also naturally acts on $\text{NS}_W(X_\eta)$ using the restriction map $\text{NS}_W(X/Y)  \twoheadrightarrow \text{NS}_W(X_\eta)$. By $[M_\eta]$ we denote the restriction of $[M]\in \text{NS}_W(X/Y)$ to $\text{NS}_W(X_\eta)$. We say a Weil divisor $M$ on $X$ is movable over $Y$ if $f_*\Os_X(M)\neq 0$ and $\text{supp}\left(\text{coker}\left(f^*f_*\Os_X(M)\to \Os_X(M)\right)\right)$ has codimension $\ge 2$, here $f: X\to Y$ is the given morphism $X/Y$.
	
	\begin{theorem}
		\label{th: finiteness_movable_orbits}
		Let $(X, B)$ be a projective good klt wlc model, $Y$ its log canonical model, and $X_\eta$ the generic fiber of $X/Y$.  In addition, we suppose $X$ is $\Q$-factorial. Then for every integer $v\in \mathbb{Z}$, the number of orbits for the natural action of $\textnormal{PsAut}(X/Y, B)$ on the following set is finite:
		\begin{equation*}
			\{ \ [M_\eta]\in \textnormal{NS}_W(X_\eta): \textnormal{vol}(M_\eta) = v,\  M \textnormal{ is a  Weil divisor on $X$ which is big and movable over $Y$ } \}.
		\end{equation*}
	\end{theorem}
	\begin{proof}
		\textbf{Step 1.} We reduce the proof to the case when $B$ is a $\Q$-divisor. Suppose the theorem holds for all projective good klt wlc models with $\Q$-coefficients. Let $B = \sum r_i\Delta_i$ be a decomposition as in \cref{lemma: log_geography}. Then $\text{PsAut}(X/Y, B) = \bigcap_i \text{PsAut}(X/Y, \Delta_i)$. In fact, the index of $\text{PsAut}(X/Y, B)$ in $\text{PsAut}(X/Y, \Delta_i)$ is bounded by $\#\supp(B)!$. Hence, the finiteness of orbits for $\text{PsAut}(X/Y, \Delta_i)$ implies the finiteness of orbits for $\text{PsAut}(X/Y, B)$. From now and, we assume $B$ is an (effective) $\Q$-divisor. 
		
		\textbf{Step 2.}
		Let $\{M_\al: \al \in I\}$ be effective divisors representing all classes in the given set, indexed by some set $I$. Then for a rational number $0<\epsilon_\al \ll 1$ the log pair $(X, B + M_\al)$ is klt of log general type over $Y$. Then we can run MMP over $Y$ with scaling \cite{BCHM10} for $(X, B + M_\al)$ which terminates with a projective $\Q$-factorial klt wlc model $\varphi_\al: (X, B + \epsilon_\al M_\al) \dashrightarrow (X_\al, B_\al + \epsilon_\al D_\al)$, where $B_\al = (\varphi_\al)_*(B)$ and $D_\al = (\varphi_\al)_* M_\al$. Since $M_\al$ is movable, it follows that all steps in the MMP are flips. Hence, $\varphi_\al$ is an isomorphism in codimension~1, and $(X_\al, B_\al)\sim_{\text{bir}} (X, B)$. Since the volume is invariant under birational morphisms and small birational maps \cite{HMX18}, we have $\text{vol}(M_\al|_F) = \text{vol}(D_\al|_{F_\al}) = v$ where $F_\al$ is a general fiber of $X_\al/Y$. Let $g_\al: X_\al \to Y_\al$ be the log canonical model for $(X_\al, B_\al + D_\al)$ over $Y$. Set $A_\al: = (g_\al)_* D_\al$ and $B_\al^\text{can}:= (g_\al)_*B_\al$. Similarly,  we obtain $\text{vol}(D_\al|_{F_\al}) = \text{vol}(A_\al|_{F_\al'})$ where $F_\al'$ is a general fiber of $Y_\al/Y$. Consequently, all pairs $(Y_\al, B_\al^\text{can})$ admit a bounded polarization by \cite{Zhu26}. Hence, the set $\|\{(Y_\al, B_\al^\text{can}): \al\in I\}\|$ is finite by \cref{th: boundedness_and_finiteness}. From Shokurov's log geography \cite[Lemma 5.3]{JJZ26} it follows that the set  $\|\{(X_\al, B_\al): \al\in I\}\|$ is finite as well. Hence, we can reduce the proof to the case when all log pairs $(X_\al,B_\al)$ are log isomorphic to some pair $(X', B')$. Then by \cref{th: finiteness_big_nef} there are at most finitely many relative big and nef divisors $D_\al$ on $X'$ up to $\text{Aut}(X'/Y, B')$ because their relative volume $\text{vol}(D_\al|_{F_\al})$ is fixed. This implies it is sufficient to prove the theorem when all $D_\al \equiv_Y D'$ for some $D'$. Clearly, for $\theta_{\al\beta} := \varphi_\beta^{-1}\circ \varphi_\al \in\text{PsAut}(X/Y, B)$ we have $(\theta_{\alpha\beta})_*M_\al \equiv_Y M_\beta$ for all $\al, \beta\in I$. This completes the proof.
	\end{proof}
	\begin{remark}
		Note that \cref{th: finiteness_big_nef} and \cref{th: finiteness_movable_orbits} stronger than the analogues statement about finiteness of orbits under the natural action of $\text{Aut}(X_\eta, B_\eta)$ and $\text{PsAut}(X_\eta, B_\eta)$ on the given sets.
	\end{remark}
	We finish the paper by showing that
		both \cref{th: finiteness_big_nef} and \cref{th: finiteness_movable_orbits} could not be strengthened to finiteness of orbits for the natural action of $\text{Aut}(X/Y, B)$ and $\text{PsAut}(X/Y, B)$ on the appropriate sets in $\text{NS}(X/Y)$ unless one fixes some global numerical data.

	\begin{example}
		Let $X$ be a smooth projective minimal surface of general type, i.e. $K_X$ is big and nef, $g: X\to Y$ be its canonical model, and $X_\eta$ be the generic fiber. Suppose $g$ contracts only one $(-2)$-curve $C$. Then for every positive integer $n\in \Zi_{>0}$ the Cartier divisor $A_n:= K_X - nC\equiv_Y -nC$ is ample over $Y$, and $\text{vol}(A_\eta) = 1$. Therefore, the number of orbits of $\text{Aut}(X/Y)$ on the following set is infinite:
		\begin{equation*}
			\{ \ [A]\in \textnormal{NS}(X/Y): \textnormal{vol}(A_\eta) = 1,\  A \textnormal{ is a Cartier divisor on $X$ which is ample over $Y$ } \}.
		\end{equation*}
		Now, let $E$ be a smooth projective genus 1 curve, and set $X_r: = X\times_k E^{\times r}$. Then one can show that the number of orbits for $\text{Aut}(X_r/Y)$ on the set
		\begin{equation*}
			\{ \ [A]\in \textnormal{NS}(X_r/Y): \textnormal{vol}(A_\eta) = 1,\  A \textnormal{ is a Cartier divisor on $X$ which is ample over $Y$ } \}.
		\end{equation*}
		is infinite for every $r\in \Zi_{\ge 0}$. By construction, $\kappa(X_r) = 2\le \dim X_r = 2+r$.
	\end{example}


\newpage
	\begin{bibdiv}
	\begin{biblist}[\small]
		
	\bib{EGAIV}{article}{
		label={EGA\,IV}
		author={Grothendieck, A.},
		author={Dieudonn\'e, J. A.},
		title={\'El\'ements de g\'eom\'etrie alg\'ebrique. {IV}. \'Etude locale des sch\'emas et des morphismes de sch\'emas. {I}},
		journal={Publ. Math. Inst. Hautes \'Etudes Sci.},
		volume={20},
		date={1964},
		pages={5--259},
	}
		
	\bib{Amb05}{article}{
		author={Ambro, F.},
		title={The moduli $b$-divisor of an lc-trivial fibration},
		journal={Compos. Math.},
		volume={141},
		date={2005},
		pages={385--403},
		number={2},
	}
	
	\bib{BBG16}{article}{
		author={Bogomolov, Fedor A.},
		author={B\"ohning, Christian},
		author={Graf von Bothmer, Hans-Christian},
		title={Birationally isotrivial fiber spaces},
		journal={Eur. J. Math.},
		volume={2},
		number={1},
		date={2016},
		pages={45--54},
	}
	
	
	\bib{BBT23}{article}{
		author={Bakker, B.},
		author={Brunebarbe, Y.},
		author={Tsimerman, J.},
		title={$o$-minimal {GAGA} and a conjecture of {G}riffiths},
		journal={Invent. Math.},
		volume={232},
		date={2023},
		pages={163--228},
		number={1},
	}
	
	\bib{Bir11}{article}{
		author={Birkar, C.},
		title={On existence of log minimal models II},
		journal={J. Reine Angew. Math.},
		volume={658},
		date={2011},
		pages={99--113},
	}
	
	\bib{Bir19}{article}{
		author={Birkar, C.},
		title={Anti-pluricanonical systems on {F}ano varieties},
		journal={Ann. of Math. (2)},
		volume={190},
		date={2019},
		pages={345--463},
		number={2},
	}
	
	
	
	\bib{Bir23}{article}{
		author={Birkar, Caucher},
		title={Geometry of polarised varieties},
		journal={Publ. Math. IH\'ES},
		volume={137},
		date={2023},
		pages={47--105},
	}
	
	\bib{BCHM10}{article}{
		author={Birkar, C.},
		author={Cascini, P.},
		author={Hacon, C. D.},
		author={McKernan, J.},
		title={Existence of minimal models for varieties of log general type},
		journal={J. Amer. Math. Soc.},
		volume={23},
		date={2010},
		pages={405--468},
		number={2},
	}
	
	\bib{BDCS24}{article}{
		author={Birkar, Caucher},
		author={Di Cerbo, Gabriele},
		author={Svaldi, Roberto},
		title={Boundedness of elliptic Calabi--Yau varieties with a rational section},
		journal={J. Differential Geom.},
		volume={128},
		number={2},
		date={2024},
		pages={463--519},
	}
	
	
	
	\bib{Bri22}{article}{
		author={Brion, Michel},
		title={Homomorphisms of algebraic groups: representability and rigidity},
		journal={Michigan Math. J.},
		volume={72},
		date={2022},
		pages={51--76},
	}
	
	\bib{CS11}{article}{
		author={Shokurov, Vyacheslav V.},
		author={Choi, Sung Rak},
		title={Geography of log models: theory and applications},
		journal={Cent. Eur. J. Math.},
		volume={9},
		date={2011},
		pages={489--534},
		number={3},
	}
	
	
	\bib{Deb01}{book}{
		author={Debarre, Olivier},
		title={Higher-Dimensional Algebraic Geometry},
		publisher={Springer},
		date={2001},
	}
	
	\bib{Fil24}{article}{
		author={Filipazzi, Stefano},
		title={On the boundedness of $n$-folds with $\kappa(X)=n-1$},
		journal={Algebr. Geom.},
		volume={11},
		number={3},
		date={2024},
		pages={318--345},
	}
	
	\bib{FHS25}{article}{
		author={Filipazzi, Stefano},
		author={Hacon, Christopher D.},
		author={Svaldi, Roberto},
		title={Boundedness of elliptic Calabi--Yau threefolds},
		journal={J. Eur. Math. Soc. (JEMS)},
		volume={27},
		number={9},
		date={2025},
		pages={3583--3650},
	}
	
	\bib{FS25}{article}{
		author={Filipazzi, Stefano},
		author={Spicer, Calum},
		title={On semi-ampleness of the moduli part},
		journal={Moduli},
		volume={2},
		date={2025},
		pages={e12, 1--27},
	}
	
	\bib{FGA05}{book}{
		editor={Fantechi, Barbara},
		editor={G\"ottsche, Lothar},
		editor={Illusie, Luc},
		editor={Kleiman, Steven L.},
		editor={Nitsure, Nitin},
		editor={Vistoli, Angelo},
		title={Fundamental Algebraic Geometry: Grothendieck's {FGA} Explained},
		publisher={American Mathematical Society},
		date={2005},
	}
	
	\bib{Fuj09}{article}{
		author={Fujino, Osamu},
		title={Effective base point free theorem for log canonical pairs---Koll\'ar type theorem},
		journal={Tohoku Math. J. (2)},
		volume={61},
		number={4},
		date={2009},
		pages={475--481},
	}
	
	
	\bib{GLSW26}{article}{
		author={Gachet, C\'ecile},
		author={Lin, Hsueh-Yung},
		author={Stenger, Isabel},
		author={Wang, Long},
		title={The effective cone conjecture for Calabi--Yau pairs},
		journal={J. Reine Angew. Math.},
		date={2026},
	}
	
	
	\bib{GW20}{book}{
		author={G\"ortz, Ulrich},
		author={Wedhorn, Torsten},
		title={Algebraic Geometry I: Schemes},
		publisher={Springer Spektrum},
		date={2020},
	}
	
	
	\bib{Har77}{book}{
		author={Hartshorne, R.},
		title={Algebraic geometry},
		publisher={Springer-Verlag},
		date={1977},
	}
	
	\bib{HH25}{article}{
		author={Hashizume, Kenta},
		author={Hattori, Masafumi},
		title={On boundedness and moduli spaces of $K$-stable Calabi--Yau fibrations over curves},
		journal={Geom. Topol.},
		volume={29},
		number={3},
		date={2025},
		pages={1619--1691},
	}
	
	\bib{HMX18}{article}{
		author={Hacon, Christopher D.},
		author={McKernan, James},
		author={Xu, Chenyang},
		title={Boundedness of moduli of varieties of general type},
		journal={J. Eur. Math. Soc. (JEMS)},
		volume={20},
		number={4},
		date={2018},
		pages={865--901},
	}
	
	
	\bib{HX15}{article}{
		author={Hacon, C. D.},
		author={Xu, C.},
		title={Boundedness of log Calabi--Yau pairs of Fano type},
		journal={Math. Res. Lett.},
		volume={22},
		date={2015},
		pages={1699--1716},
		number={6},
	}
	
	
	\bib{HPT18}{article}{
		author={Hassett, Brendan},
		author={Pirutka, Alena},
		author={Tschinkel, Yuri},
		title={Stable rationality of quadric surface bundles over surfaces},
		journal={Acta Math.},
		volume={220},
		number={2},
		date={2018},
		pages={341--365},
	}
	
	\bib{Huy18}{article}{
		author={Huybrechts, Daniel},
		title={Finiteness of polarized {K3} surfaces and hyperk\"ahler manifolds},
		journal={Ann. H. Lebesgue},
		volume={1},
		date={2018},
		pages={227--248},
	}
	
	
	
	\bib{Jia25}{article}{
		author={Jiao, Junpeng},
		title={Boundedness of polarized log Calabi--Yau fibrations},
		journal={J. Differential Geom.},
		volume={130},
		number={3},
		date={2025},
		pages={635--675},
	}
	
	\bib{JJZ26}{article}{
		author={Jiang, Xiaowei},
		author={Jiao, Junpeng},
		author={Zhu, Minzhe},
		title={Boundedness of polarized log Calabi--Yau fibrations with bounded bases},
		journal={J. Eur. Math. Soc.},
		date={2026},
	}
	
	\bib{Kaw78}{article}{
		author={Kawamata, Yujiro},
		title={On deformations of compactifiable complex manifolds},
		journal={Math. Ann.},
		volume={235},
		date={1978},
		pages={247--265},
		number={3},
	}
	
	\bib{Kaw85}{article}{
		author={Kawamata, Y.},
		title={Minimal models and the Kodaira dimension of algebraic fiber spaces},
		journal={J. Reine Angew. Math.},
		volume={363},
		date={1985},
		pages={1--46},
	}
	
	\bib{Kaw97}{article}{
		author={Kawamata, Yujiro},
		title={On the cone of divisors of Calabi--Yau fiber spaces},
		journal={Int. J. Math.},
		volume={8},
		number={5},
		date={1997},
		pages={665--687},
	}
	
	
	\bib{Kaw08}{article}{
		author={Kawamata, Yujiro},
		title={Flops connect minimal models},
		journal={Publ. Res. Inst. Math. Sci.},
		volume={44},
		date={2008},
		pages={419--423},
		number={2},
	}
	
	\bib{KM98}{book}{
		author={Koll\'{a}r, J.},
		author={Mori, S.},
		title={Birational Geometry of Algebraic Varieties},
		publisher={Cambridge University Press},
		date={1998},
	}
	
	\bib{Kol93}{article}{
		author={Koll\'ar, J\'anos},
		title={Effective base point freeness},
		journal={Math. Ann.},
		volume={296},
		date={1993},
		pages={595--605},
		number={4},
	}
	
	\bib{Kol96}{book}{
		author={Koll\'ar, J.},
		title={Rational Curves on Algebraic Varieties},
		publisher={Springer-Verlag},
		date={1996},
	}
	
	\bib{Kol13}{book}{
		author={Koll\'ar, J\'anos},
		title={Singularities of the Minimal Model Program},
		publisher={Cambridge University Press},
		date={2013},
	}
	
	\bib{Kol23}{book}{
		author={Koll\'ar, J\'anos},
		title={Families of Varieties of General Type},
		publisher={Cambridge University Press},
		date={2023},
	}
	
	
	\bib{Kol26}{article}{
		author={Koll\'ar, J.},
		title={N\'eron models, minimal models, and birational group actions},
		journal={Pure Appl. Math. Q.},
		volume={22},
		date={2026},
		pages={647--673},
		number={2},
	}
	

	
	\bib{Lan24}{article}{
		author={Landesman, Aaron},
		title={Invariance of the tame fundamental group under base change between algebraically closed fields},
		journal={Essent. Number Theory},
		volume={3},
		number={1},
		date={2024},
		pages={1--18},
	}
	
	\bib{Laz04}{book}{
		author={Lazarsfeld, Robert},
		title={Positivity in Algebraic Geometry {I}},
		publisher={Springer},
		date={2004},
	}
	
	\bib{Li25}{article}{
		author={Li, Jennifer},
		title={A cone conjecture for log Calabi--Yau surfaces},
		journal={Forum Math. Sigma},
		volume={13},
		date={2025},
		pages={e15},
	}
	
	
	\bib{Mat02}{book}{
		author={Matsuki, K.},
		title={Introduction to the Mori Program},
		publisher={Springer},
		date={2002},
	}
	
	\bib{Mor93}{incollection}{
		author={Morrison, David R.},
		title={Compactifications of moduli spaces inspired by mirror symmetry},
		booktitle={Journ\'ees de g\'eom\'etrie alg\'ebrique d'Orsay (Juillet 1992)},
		number={218},
		date={1993},
		pages={243--271},
		publisher={Soci\'et\'e Math\'ematique de France},
	}
	
	\bib{SS19}{book}{
		author={Sch\"utt, Matthias},
		author={Shioda, Tetsuji},
		title={Mordell--Weil lattices},
		publisher={Springer},
		date={2019},
	}
	

	
	\bib{Ser26a}{article}{
		author={Serebrennikov, Daniil},
		title={On the finiteness of log surfaces},
		journal={Int. J. Math.},
		volume={37},
		date={2026},
		number={13},
		pages={2650065}
	}
	
	\bib{Ser06}{book}{
		author={Sernesi, E.},
		title={Deformations of Algebraic Schemes},
		publisher={Springer},
		date={2006},
	}
	
	\bib{Ser58}{article}{
		author={Serre, Jean-Pierre},
		title={Espaces fibr\'es alg\'ebriques},
		journal={S\'eminaire Claude Chevalley},
		volume={3},
		number={1},
		date={1958},
		pages={1--37},
	}
	
	
	
	\bib{Ste76}{book}{
		author={Steenbrink, J. H. M.},
		title={Mixed {H}odge Structure on the Vanishing Cohomology},
		publisher={Department of Mathematics, University of Amsterdam},
		date={1976},
	}
	
	\bib{Sze99}{article}{
		author={Szendr\H{o}i, Bal\'azs},
		title={Some finiteness results for Calabi--Yau threefolds},
		journal={J. Lond. Math. Soc. (2)},
		volume={60},
		date={1999},
		pages={689--699},
		number={3},
	}
	
	
	
	\bib{Tot10}{article}{
		author={Totaro, B.},
		title={The cone conjecture for Calabi--Yau pairs in dimension 2},
		journal={Duke Math. J.},
		volume={154},
		date={2010},
		pages={241--263},
		number={2},
	}
	
	\bib{Zhu25}{article}{
		author={Zhu, Minzhe},
		title={Boundedness of stable minimal models with klt singularities},
		journal={Int. Math. Res. Not. IMRN},
		volume={2025},
		number={2},
		date={2025},
		pages={rnae293},
	}
	
	
	\bib{BFMT25}{misc}{
		author={Bakker, Benjamin},
		author={Filipazzi, Stefano},
		author={Mauri, Mirko},
		author={Tsimerman, Jacob},
		title={Baily--Borel compactifications of period images and the b-semiampleness conjecture},
		date={2025},
		note={ArXiv e-print, \href{https://arxiv.org/abs/2508.19215}{arXiv:2508.19215}.},
	}
	
	
	
	\bib{Bri18}{misc}{
		author={Brion, Michel},
		title={Notes on automorphism groups of projective varieties},
		date={2018},
		note={\href{https://www-fourier.univ-grenoble-alpes.fr/~mbrion/autos_final.pdf}{Lecture notes}.},
	}
	
	

	
	\bib{Ser26b}{misc}{
		author={Serebrennikov, Daniil},
		title={Constructibility aspects of the cone conjecture},
		date={2026},
		note={ArXiv e-print, \href{https://arxiv.org/abs/2604.27303}{arXiv:2604.27303}.},
	}
	
	\bib{Ser26c}{misc}{
		author={Serebrennikov, Daniil},
		title={Towards the global Torelli theorem},
		date={2026},
		note={ArXiv e-print, \href{https://arxiv.org/abs/2608.27441}{arXiv:2608.27441}.},
	}
	
	\bib{Zhu26}{misc}{
		author={Zhu, Haoran},
		title={Boundedness of polarized log Calabi--Yau fibrations},
		date={2026},
		note={ArXiv e-print, \href{https://arxiv.org/abs/2610.13178}{arXiv:2610.13178}.},
	}
	
	\bib{Stacks}{misc}{
		label={Stacks},
		author={{Stacks Project Authors}},
		title={The Stacks Project},
		date={2026},
		note={\url{https://stacks.math.columbia.edu}},
	}

	\end{biblist}
\end{bibdiv}
\end{document}